\documentclass[12pt]{article}

\usepackage{amsmath, amsthm, amssymb, amsfonts, enumerate, dsfont}
\usepackage[applemac]{inputenc} % accents 
\usepackage[english]{babel}
\usepackage{tikz}

\usetikzlibrary{arrows,automata,shapes}

\def \R{I\!\!R}
\def \E{I\!\!E}

\def \N{I\!\!N}

\newtheorem{thm}{Theorem}[section]
\newtheorem{cor}[thm]{Corollary}
\newtheorem{lem}[thm]{Lemma}
\newtheorem{pro}[thm]{Proposition}
\newtheorem{defi}[thm]{Definition}

\newtheorem{exa}[thm]{Example}

\def \1{\mathbf{1}}

\newcommand{\m}{\mathbb}
\def\abstract{\begin{center} \small\bf Abstract\end{center}\small}

 \title{Hidden Stochastic Games    
  and Limit Equilibrium Payoffs}

 \author{J\'er\^ome Renault$^*$, Bruno Ziliotto\thanks{TSE (GREMAQ, Universit\' e Toulouse 1 Capitole), 21 all\' ee de Brienne, 31000
Toulouse, France. Both authors gratefully  acknowledge   the support of the Agence Nationale de la Recherche, under grant ANR JEUDY, ANR-10-BLAN 0112, and thank J. Bolte, T. Mariotti  and T. Tomala for fruitful discussions.}
}
\date{\today}

\begin{document}     

\bibliographystyle{ieeetr}
\maketitle

\begin{abstract}
 We consider 2-player stochastic games with perfectly observed actions, and study the limit, as the discount factor goes to one, of the  equilibrium payoffs set. In the usual setup where current states are observed by the players,  we  first  show that the set of {\it stationary} equilibrium payoffs always converges.  We then  provide the first examples   where  the whole set  of equilibrium payoffs diverges. The construction can be  robust to   perturbations of the payoffs, and to the introduction of normal-form correlation. 
 Next we naturally  introduce the   more general model of hidden stochastic game,  where the players publicly receive imperfect signals   over current states. In this setup  we present a last  example  where not only  the limit set  of equilibrium payoffs does not exist, but there is no converging selection of equilibrium payoffs. The   example is  symmetric and   robust in many aspects, in particular   to the introduction of extensive-form correlation  or communication  devices. No uniform equilibrium payoff exists, and the equilibrium set has full support for each discount factor and each initial state.   \end{abstract}

\section{Introduction}
 
 Most economic and social interactions have a dynamic aspect, and equilibrium plays of dynamic games are typically not obtained by successions of  myopic equilibria of the current one-shot interaction, but need to take into account  both  the effects of actions over current payoffs, and over future payoffs in the continuation game. In this paper we consider     dynamic games with 2 players\footnote{This simplifies the exposition, but  our results extend to the $n$-player case.}, where the actions taken by the players are perfectly observed at the end of every stage. We  denote by  $E_\delta$, resp. $E'_\delta$,  the set  of Nash equilibrium payoffs, resp. sequential equilibrium payoffs,  of the $\delta$-discounted game, and we write  $E_\infty$ for the  set of uniform equilibrium payoffs of the dynamic game. We  mainly study the limit\footnote{for the Hausdorff distance, defined as $d(A,B)= \max\{\max_{a \in A} d(a,B),\max_{b \in B}d(b,A) \} $ for   $A$ and $B$ non-empty compact subsets of $\R^2$.  $d(A,B)\leq \varepsilon$ means that:  every point in $A$ is at a distance at most $\varepsilon$ from a point in $B$, and conversely.}  of $E'_\delta$ as players get extremely patient, i.e. as the discount factor goes to one.  
 
In a repeated game,  the dynamic interaction consists of the repetition of a given one-shot game, and we have the standard Folk Theorems, with pioneering work from the seventies  by Aumann and Shapley, and Rubinstein.   Regarding sequential equilibrium payoffs, the Folk theorem of Fudenberg and Maskin (1986) implies  that for generic payoff functions, $E'_\delta$  converges to the set of feasible and individually rational payoffs of the one-shot game, and Wen (1994) showed how to adapt the notion of individually rational payoffs to obtain a Folk theorem without genericity assumption. Without assumptions on the payoffs, $E_\infty$ coincides with the set of   feasible and individually rational payoffs of the one-shot game, and the set of Nash  equilibrium payoffs  $E_\delta$ also converges to this set (see Sorin (1986)). These results have been  generalized in many ways to games with imperfectly observed actions  (see e.g. Abreu {\it et al.} (1990), Fudenberg Levine (1991), Fudenberg Levine Maskin (1994), Fudenberg {\it et al.} (2007), Lehrer (1990, 1992a, 1992b), or Renault Tomala (2004, 2011)), but this is beyond the scope of the present paper.

Stochastic games were introduced by Shapley (1953) and generalize repeated games: the   payoff functions of the players evolve from stage to stage, and depend  on a state variable observed by the players,  whose evolution is influenced by the players' actions.  In the zero-sum case,  Bewley and Kohlberg (1976) proved the existence of the limit of  the discounted value (hence of $E_\delta$ and $E'_\delta$)  when $\delta $ goes to one.  An example of Sorin (1984) shows  that in the general-sum case $\lim_{\delta\to 1} E_\delta$ and $E_\infty$ may be non-empty and disjoint. Vieille (2000) proved that $E_{\infty}$ is always non-empty, that is there   exists a uniform equilibrium payoff\footnote{The   generalization of this result to more players is  a well-known    open question in dynamic games.}.  

Regarding discounted equilibrium payoffs in stochastic games, several Folk theorems  have been  proved under various assumptions. Dutta (1995) assumes that the set of  long-run feasible payoffs is independent of the initial state, has full dimension, and that minmax long-run payoffs also do not depend on the initial state. Fudenberg\footnote{Fudenberg and Yamamoto 2011, as well as H\"{o}rner {\it et al.} (2011), consider the more general case of imperfect public monitoring.}  and Yamamoto (2011)  assume that the stochastic game is irreducible (all players but one can always drive the current state where they want, possibly in many stages, with positive probability). H\"{o}rner {\it et al.} (2011) generalize the recursive methods of Fudenberg Levine Maskin (1994) to compute a limit equilibrium set in stochastic games with imperfect public monitoring when this limit set does not depend on the initial state (this happens when the Markov chain induced by any Markov strategy profile is irreducible). 

All the above assumptions somehow require  that the stochastic game does not depend too much on the initial state, and in particular forbid the existence of multiple absorbing\footnote{When an absorbing state is reached, the play will remain forever in that state,  no matter the actions played.}  states with different equilibrium payoffs. We believe that it is also meaningful to study stochastic games where the actions taken can have irreversible effects on future plays. This is    the case in many situations, for example in  stopping games when each player only acts once and has to decide when to do so, or when the actions partially represent investment decisions, or extractions of exhaustible resources. 

The first contributions of our paper concern  2-player stochastic games with finitely many states and actions. We first prove that  the set of {\it stationary}\footnote{A stationary strategy of a player plays after every history  a mixed action which only depends on the current state.} equilibrium payoffs of the discounted game always converges to a non-empty set.  As a consequence,  there always exists a  selection   of  $E'_\delta$ which converges. Secondly, we show that the convergence property  cannot be extended to  the Nash or sequential equilibrium payoffs sets by providing the first examples of    stochastic games  where neither   $E_\delta$ nor  $E'_\delta$  converges: the limit of the equilibrium set may simply   not exist in a stochastic game. Our first example is robust to the introduction of normal-form correlation, and our second example shows that the non convergence property can be robust to  small perturbations of the payoffs. However we point out that both examples are not robust to  the introduction of an extensive-form correlation device. In each case, the set of equilibrium payoffs has empty interior for any  discount factor, and the limits of converging selections of  $E'_\delta$ coincide with the uniform equilibrium payoffs. In these examples, we believe that the elements of the limit set of stationary equilibrium payoffs  emerge as   natural outcomes of the game played by extremely patient players.   \\

In the rest of the paper we introduce the more general model of {\it hidden} stochastic games, and we refer to the above original model as {\it standard} stochastic games. In a hidden stochastic game, the players still perfectly observe past actions but no longer perfectly observe current states, and rather  receive at the beginning of every stage a public, possibly random, signal on the current state. So players have incomplete information over the sequence of states, but this information is common to both players. Note that in the zero-sum case this model was considered in Gimbert et al. $\cite{G14}$. Hidden stochastic games are generalizations of hidden Markov decision processes (where there is a single agent), hence the name. Hidden stochastic games also generalize repeated games with common incomplete information on the state. We believe this model is meaningful in many interactions where the fundamentals are not perfectly known to the players. We present in particular two examples of economic interactions that could  be modeled as a hidden stochastic game. The first example  is a Cournot competition on a market for a natural exhaustible resource, where the players have common incomplete information on the stock of natural resources remaining. The second example is an oligopoly competition with a single good (e.g., traditional  tv sets)  where the state variable includes  parameters known to the firms, such as the current demand level, but also parameters imperfectly known such as the trend of the market or the overall state of the economy. %or the  possible introduction of close other goods (e.g., flat screen tv) which can potentially have a dramatic impact on future demands for the good, but is difficult to predict in advance by the firms. 

Surprisingly enough,  few papers have already considered stochastic games with imperfect observation of the state. In the zero-sum context, Venel (2012) studied  hidden stochastic games where the players do not receive any signal   on the state  during the game, and proves under a commutativity assumption over transitions the existence of the limit value, as well as the stronger notion of uniform value (corresponding to uniform equilibrium payoffs). Ziliotto (2013) showed that the commutativity assumption was needed for  Venel's result, and provided an example of a zero-sum hidden stochastic game\footnote{The example of Ziliotto provides a negative answer to two  conjectures of Mertens (1986) for zero-sum dynamic games.}
 with payoffs in $[0,1]$ where the $\delta$-discounted value oscillates between $1/2$ and $5/9$ when $\delta$ goes to one.   \\
 
Given   parameters $\varepsilon$ in $(0,5/12)$ and $r$ in $(0,\varepsilon/5)$, we provide an example of a 2-player hidden stochastic game with all payoffs in $[0,1]$, four actions for each player,   having    the following features:

- the game is symmetric between the players, 

- the players have incomplete information over the current state, but the public signals received are informative enough for the players to know the current stage payoff functions at the beginning of every stage. As a consequence, the players know their current payoffs during the play. 

- there are 13 states, and for any initial state and discount factor the set of sequential equilibrium payoffs contains a square with side $2r$, hence  has full dimension.

- for a specific initial state $k_1$, there exist subsets $\Delta_1$ and $\Delta_2$ of discount factors, both containing 1 as a limit point, such that
for all discount factors in $\Delta_1$,   the corresponding set of sequential equilibrium payoffs is exactly the square  $E_1$ centered in $(\epsilon,\epsilon)$ with side $2r$, whereas 
for all discount factors in $\Delta_2$, the set\footnote{As an illustration, if $\varepsilon=.3$ and $r=.05$, for any discount factor in $\Delta_1$ the set of equilibrium payoffs is the square  $E_1={[.25,.35]}^2$, and for any discount in $\Delta_2$ the set of equilibrium payoffs is the square   $E_2={[.65,.75]}^2$.} of sequential equilibrium payoffs is the square $E_2$ centered in $(1-\epsilon,1-\epsilon)$ with side $2 r$. In each case the associated square  is also the set of Nash equilibrium payoffs, the set of (normal or extensive-form) correlated equilibrium payoffs, and the set of communication\footnote{ introduced in Myerson, 1986 and Forges, 1986.}  equilibrium payoffs  of the discounted game. Since these two squares are disjoint, there is no converging  selection of equilibrium payoffs, and the game has no uniform equilibrium payoff.

Moreover the example is robust to small perturbations of the payoffs: if one perturbs all payoffs of the game by at most $\frac{1}{2} r (\varepsilon-5r)$, the set of discounted equilibrium payoffs of the perturbed game with initial state $k_1$ still does not converge,  no  converging selection of equilibrium payoffs exists and there is no uniform equilibrium payoff. \\

 Our last  example is thus robust in many aspects, and it seems impossible to affect to this game a reasonable  limit equilibrium payoff. The model of hidden stochastic games may be seen as a small departure from the standard model of stochastic game, but it seems very difficult for an expert to find any good answer to the informal question: ``The game being played by extremely patient players, which outcome is likely to form ?"\\

 We study  standard stochastic games in section \ref{secSG}. Hidden stochastic games are introduced in section \ref{secHSG}, and our last example is presented in section \ref{secProof}. The construction   elaborates  and  improves on the zero-sum construction of Ziliotto (2013). The  presentation is done here in 5 progressive steps,  starting with a Markov chain on $[0,1]$, then a Markov Decision Process, then a zero-sum stochastic game with infinite state space,  a zero-sum hidden stochastic game and a final  example. A few proofs are relegated to   the Appendix.  \\

 We  denote  respectively by  $\N$, $\R$  and $\R_+$ the sets of non negative integers, real numbers and non negative real numbers.

\section{Standard Stochastic Games}\label{secSG}

We consider a 2-player stochastic game. Let $K$, $I$ and $J$ respectively be the finite sets of states, actions for player 1 and actions for player 2. $k_1$ in $K$ is the initial state, $u_1$ and $u_2$ are the state  dependent utility functions from $K\times I \times J$ to $\R$, and $q$   is the transition function from $K \times I \times J$ to $\Delta(K)$,   the set of probabilities over  $K$. At every period $t \geq 1$ players first learn the current state $k_t\in K$ and  simultaneously select actions $i_t\in I$ and $j_t\in J$. These actions are then publicly observed, the stage payoffs are $u_1(k_t, i_t,j_t)$ for player 1 and $u_2(k_t, i_t,j_t)$ for player 2, a new state $k_{t+1}$ is selected according to the distribution $q(k_t,i_t,j_t)$, and the play goes to the next period. Given a discount factor $\delta$ in $[0,1)$, the $\delta$-discounted stochastic game  is the infinite horizon  game where   player 1 and player 2's payoffs  are respectively  $(1-\delta) \sum_{t=1}^\infty \delta^{t-1} u_1(k_t,i_t,j_t)$ and $(1-\delta) \sum_{t=1}^\infty \delta^{t-1} u_2(k_t,i_t,j_t)$.

Let $E_\delta$ and $E'_\delta$ respectively denote the set of Nash equilibrium payoffs and the   set of perfect\footnote{subgame-perfect equilibrium, or equivalently here, sequential  equilibrium, or equivalently  perfect public equilibrium, as defined in Fudenberg Levine Maskin 1994 for repeated games with perfect public monitoring and extended to stochastic games in Fudenberg Yamamoto 2011.} equilibrium payoffs of the  $\delta$-discounted stochastic game.  Standard fixed-point   arguments   show the existence of a stationary  equilibrium in this game, and the associated equilibrium payoff lies in $E'_\delta$ and $E_\delta$. 
In this paper, we  are mainly interested in the asymptotic behavior  of these  sets when players become more and more patient, i.e. we will look for their limit\footnote{All limits of sets in the paper are to be understood for the Hausdorff distance between non-empty compact sets of $\R^2$.}    when $\delta$ goes to 1. And we  will also  briefly consider the set $E_\infty$ of uniform equilibrium\footnote{Throughout the paper, we say that a vector $u$  in $\R^2$ is a uniform equilibrium payoff  if   for all $\varepsilon>0$, there exists a strategy profile such that for all high enough  discount factors, the profile is a $\varepsilon$-Nash equilibrium of the discounted game with payoff $\varepsilon$-close to $u$, see Sorin 1986, Mertens Sorin and Zamir 1994 or Vieille 2000 for related definitions.}  payoffs of the stochastic game.\\

When there is a single state, the game is a standard repeated game with perfect monitoring, and we have  well-known    Folk Theorems.

For   zero-sum stochastic games, Shapley (1953) proved that the value $v_\delta$ exists and players have stationary optimal strategies, so $E_\delta$ and $E'_\delta$ are singletons. Bewley and Kohlberg (1976)  proved the convergence of $v_\delta$ (hence, of $E_\delta$ and $E'_\delta$) using algebraic arguments\footnote{Recently, Bolte {\it et al.} (2013) extended this  algebraic approach to a class of  stochastic games with infinite action sets.}.

The following proposition shows how the Bewley Kohlberg result extends  to general-sum games. Denote by  $E''_\delta$ the set  of {\it stationary} equilibrium payoffs of the $\delta$-discounted game.  In the zero-sum case, $E_\delta=E'_\delta=E''_\delta$. In general,    stationary equilibria are very simple equilibria where the strategies of the players are particularly restricted, and  $E''_\delta$ is a  subset of $E'_\delta$.  

\begin{pro} \label{pro1} There exists a non-empty compact set $E$ such that: $$E''_{\delta} \xrightarrow[\delta \to 1]{} E.$$ 
\end{pro}

\noindent In the case of repeated games (a single state),  $E$ reduces to the  set of mixed Nash equilibrium payoffs of the one-shot game, hence may not be convex. The proof of proposition \ref{pro1} is in the Appendix and largely relies on the semi-algebraicity of the set of discount factors and associated stationary equilibria and payoffs. As stated here, it holds for  any 2-player stochastic game with finitely many states and actions, but the proof easily extends to the $n$-player  case. As a consequence, using a point in $E$ one can construct  a selection  of $(E'_\delta)_\delta$ which converges, i.e. it is possible to select, for each discount $\delta$, a perfect  equilibrium payoff $x_\delta$ of the corresponding game in a way such that $x_\delta$ has a limit when $\delta$ goes to one. 
\begin{cor}
There exists a  converging  selection of $(E'_\delta)_\delta$.
\end{cor}

\noindent This corollary can also be easily deduced from Mertens Sorin Zamir (1994, Lemma 2.2 in chapter 7) and Neyman (2003, Theorem 5), who have proved the existence of a semi-algebraic selection of $E''_{\delta}$. Because payoffs are bounded, this selection converges.
%We mention that  this  existence of a converging selection of equilibrium payoffs  can also be  easily deduced from  . 
\\

It is then now natural to ask if the convergence property of proposition \ref{pro1} also holds for $E_\delta$ and $E'_\delta$. We conclude this section by providing the first  examples of  stochastic games  where these  sets of  equilibrium payoffs  diverge.

\begin{pro} \label{pro2}
There exists a 2-player stochastic game where neither $E_\delta$ nor $E'_\delta$ converge. The construction can be  robust to   perturbations of the payoffs, and to the introduction of normal-form correlation. 
\end{pro}

To prove the proposition, we first present  a simple example (Example \ref{exa1} below) where $E_\delta$ and  $E'_\delta$ diverge, and which is robust to the introduction of normal-form correlation.  Then we provide a more elaborate example (Example \ref{exa1.5}) which is robust to    perturbations of the payoffs. \\

\begin{exa}\label{exa1}\rm  Consider the  stochastic game   represented by the following picture.
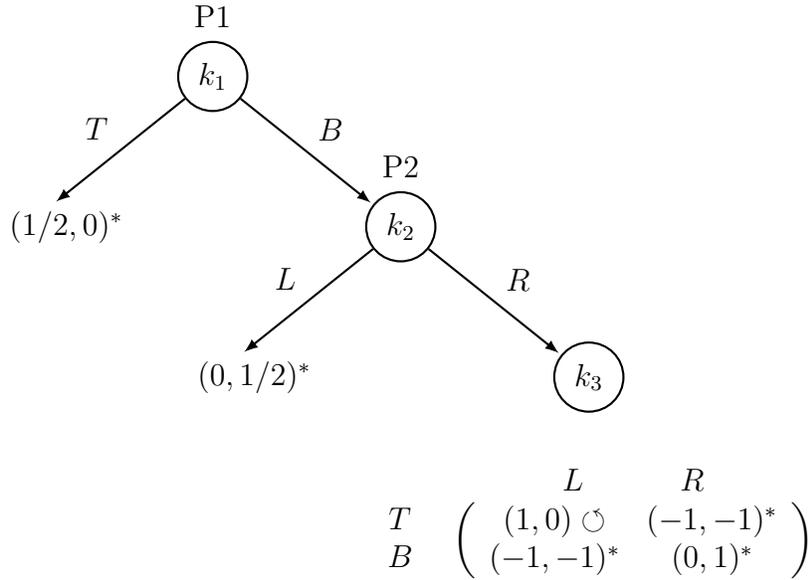
\unitlength 0,7mm
% sp�cifie l'unit�

\begin{center}
\begin{tikzpicture}[>=stealth',shorten >=1pt,auto,node distance=4cm,thick,main
 node/.style={circle,draw,font=\Large\bfseries}]

\node [draw,text width=0.4cm,text centered,circle] (O) at (0,0) {$k_1$};
\node [draw,text width=0.4cm,text centered,circle] (A) at (2.5,-2) {$k_2$};
\node [draw,text width=0.4cm,text centered,circle] (B) at (5,-4) {$k_3$};
\node [text width=0.4cm,text centered,circle] (C) at (-2.5,-2) {$(1/2,0)^*$};
\node [text width=0.4cm,text centered,circle] (D) at (0,-4) {$(0,1/2)^*$};

\draw[->,>=latex] (O) to[] node[midway, above left] {$T$}(C);
\draw[->,>=latex] (O) to[] node[midway,above right] {$B$}(A);
\draw[->,>=latex] (A) edge[] node[midway,above right] {$R$}(B);
\draw[->,>=latex] (A) edge[] node[midway, above left] {$L$}(D);

\draw(O)+(0,0.5) node[above]{P1};
\draw(A)+(0,0.5) node[above]{P2};
\draw (B)+(0,-1) node[below]{
\vspace{4cm}
$\begin{array}{cc}
 \; & \begin{array}{cc}
 L\;\; \;\;&\; \; \; \; R \\
 \end{array} \\
 \begin{array}{c}
 T \\
 B \\
 \end{array} &
 \left( \begin{array}{cc}
  (1,0) \circlearrowleft& (-1,-1)^* \\
  (-1,-1)^* & (0,1)^* \\
\end{array}\right)\\
\end{array} $};
%\draw[->,>=latex] (A) to[bend right] node[midway, left] {$W_1,1-\alpha,s_1$} (O);

%\draw[->,>=latex] (O) to[] node[midway, right] {$J_1,s_0^*$}(B);

\end{tikzpicture}
\end{center}

\vspace{1cm}

\noindent There are 7 states: $k_1$ (the initial state), $k_2$, $k_3$ and 4 absorbing states: $(1/2,0)^*$, $(0,1/2)^*$, $(-1,-1)^*$ and $(0,1)^*$. When an absorbing state $(a,b)^*$ is reached, the game stays there forever and at each stage the payoffs to player 1 and player 2 are respectively $a$ and $b$. The sets of actions are $I=\{T,B\}$ for player 1 and $J=\{L,R\}$ for player 2. The transition from state $k_1$
 only depends on player 1's action, as indicated in the above figure, and similarly the transition from state $k_2$
 only depends on player 2's action.  If in state $k_3$ the action profile $(T,L)$ is played, the vector payoff is $(1,0)$ and the play remains in $k_3$. To conclude the description, we have to specify the   payoffs in states $k_1$, $k_2$, and $k_3$. The payoff in    $k_1$ is $(1/2,0)$ if $T$ is played and $(1/2,1/2)$ if $B$ is played. The payoff  in state  $k_2$   does not depend on the actions played and   is    $(1/2,1/2)$, and the payoffs in state $k_3$ are simply given by the bimatrix $ 
 \left( \begin{array}{cc}
  (1,0)    & (-1,-1)  \\
  (-1,-1)  & (0,1)  \\
\end{array}  \right)$. 
 
 For each discount, it is clear that $(1/2,0)$ is in $E_\delta$, and the question is whether there are other equilibrium payoffs, for instance $(1/2,1/2)$.
 
First consider any $\delta$ in $[0,1)$, and  a Nash equilibrium $(\sigma, \tau)$ of the $\delta$-discounted stochastic game    with equilibrium payoff $(x,y)$.  Because Player 1 can play $T$ in the initial state, $x\geq 1/2$. Because the sum of payoffs never exceeds 1, we have $x+y\leq 1$. Assume now that under $(\sigma, \tau)$ the state $k_3$ has positive probability to be reached, and denote by $(x_3,y_3)$ the  discounted payoffs  induced by $(\sigma, \tau)$ given that   $k_3$ is reached. We have $x_3\geq 1/2$, because player 1 will  not accept to play $B$ at $k_1$ if he obtains a payoff lower than 1/2 afterwards. Similarly, $y_3\geq 1/2$. Since $x_3+y_3\leq 1$, we get $x_3=y_3=1/2$, so  $(1/2,1/2)$ is an equilibrium payoff of the reduced stochastic game: 

\centerline{  $\begin{array}{cc}
 \; & \begin{array}{cc}
 L\;\; \;\;&\; \; \; \; R \\
 \end{array} \\
 \begin{array}{c}
 T \\
 B \\
 \end{array} &
 \left( \begin{array}{cc}
  (1,0) \circlearrowleft & (-1,-1)^* \\
  (-1,-1)^* & (0,1)^* \\
\end{array}\right)\\
\end{array}$ .} 
\vspace{0.5cm}
\noindent The unique way to obtain $(1/2,1/2)$ as a feasible payoff in the reduced game is to play first $(T,L)$ for a certain  number of periods $N$, then $(B,R)$ at period $N+1$. Given $\delta$, the integer $N$ has to satisfy $(1-\delta) \sum_{t=1}^{N} \delta^{t-1}= 1/2$, that is:
$$\delta^{N}=\frac{1}{2}.$$ If no such integer $N$ exists, we obtain that $(1/2,1/2)$ is not an equilibrium payoff of the reduced game, so  under $(\sigma, \tau)$ the state $k_3$ has zero probability to be reached, which implies that $x=1/2$ and $y=0$. \\

 We define $\Delta_1$ as the set of discount factors of the form $\delta= {(\frac{1}{2})}^{1/N}$, where $N$ is a positive integer, and we put $\Delta_2=[0,1)\backslash \Delta_1$. We have obtained:\\
 
 %  Put  $\Delta_2$ be the complementary of $\Delta_1$ within $[0,1)$.
 
\begin{lem} \label{lem1} For all $\delta$ in $\Delta_2$, $E_\delta=E'_\delta=\{(1/2,0)\}.$ \end{lem}
 
 Consider now $\delta$ in $\Delta_1$, and  $N$  such that $\delta^{N}=\frac{1}{2}.$  The pure strategy profiles where: $T$ is played at stage 1,  $R$ is played at stage 2, $(T,L)$ is played for $N$ periods from stage 3 to stage $N+2$, and $(B,R)$ is played at stage $N+3$, form a  subgame-perfect Nash equilibrium of the $\delta$-discounted game with payoff $(1/2,1/2)$. By mixing between $T$ and $B$ in $k_1$, it is then possible to obtain any point $(1/2,x)$, with $0\leq x \leq 1/2$, as an equilibrium payoff. And no other point can be obtained, because  in every equilibrium, the vector payoff  conditional on  $k_3$ being  reached, is  $(1/2,1/2)$. We have obtained:\\
 
\begin{lem} \label{lem2} For all $\delta$ in $\Delta_1$, $E_\delta=E'_\delta=\{1/2\}\times [0,1/2]$. \end{lem}
 
 Because both $\Delta_1$ and $\Delta_2$ contain discount factors arbitrarily close to 1, lemmas \ref{lem1} and \ref{lem2}  establish that neither $E_\delta$ nor $E'_\delta$ converge\footnote{One can also consider for any positive integer $n$, the set of Nash equilibrium payoffs  $E_n$ and subgame-perfect equilibrium payoffs $E'_n$ of the $n$-period stochastic game, where the overall payoff is defined as the arithmetic average of the stage payoffs. Similar arguments show that in  Example \ref{exa1}, we have $E_n=E'_n=\{1/2\}\times [0,1/2]$ for   $n$ even, and  $E_n=E'_n=\{(1/2,0)\}$ for   $n$ odd. So $E_n$ and $E'_n$ also do not converge when $n$ goes to infinity.}. 
 
 Consider now normal-form correlated equilibrium payoffs, i.e. Nash equilibrium payoffs of games where the players may initially   receive   private signals independent of the payoffs. For $\delta$ in $\Delta_2$,  the proof of Lemma \ref{lem1} applies and the set of normal-form correlated equilibrium payoffs still is the singleton $\{(1/2,0)\}$. So the set of normal-form correlated equilibrium payoffs can not converge when the discount factor goes to one. 
 
 Notice that  an important feature  of  the reduced game is that there is (at most) a unique way to obtain the payoff (1/2,1/2). As soon as one perturbs the payoffs, this property will disappear, and  example \ref{exa1} is not robust to perturbations of the payoffs of the stochastic game. \hfill $\Box$ \end{exa}

\vspace{0,5cm}
 
%
%We now present a second example which is robust to  perturbations of the payoffs.  

\begin{exa} \label{exa1.5} \rm   The  stochastic game  is  represented by the following matrices:
\vspace{0,5cm}

State $k_1$: $\begin{array}{cc}
 \; & \begin{array}{ccc}
 L \hspace{1,8cm} &M\;\;\;\;\;\;\;\;\;&\;\;\;\;R\;\;\;\;\\
 \end{array} \\
 \begin{array}{c}
 T \\
 B \\
 \end{array} &
 \left( \begin{array}{ccc}
  (1,1)^* & (-30,-30)^*&(-30,-30)^* \\
  (-1,-1)^{k_2}& (-30,-30)^*& (-30,-30)^* \\
\end{array}\right)\\
\end{array}  $

\vspace{0,5cm}

State $k_2$: $\begin{array}{cc}
 \; & \begin{array}{ccc}
 L \hspace{1,8cm} &M\;\;\;\;\;\;\;\;\;&\;\;\;\;R\;\;\;\;\\
 \end{array} \\
 \begin{array}{c}
 T \\
 B \\
 \end{array} &
 \left( \begin{array}{ccc}
  (-1,-1)\circlearrowleft & (-12,-11)^*& (-4,-7)^* \\
  (-22,-12)^*& (3,-2)^*& (-9,-4)^* \\
\end{array}\right)\\
\end{array}  $

%\begin{center}
%\begin{picture}(105,70)
%\put(35,65){\line(-2,-1){70}}
%\put(35,65){\line(2,-1){70}}
% 
% 
%\put(20,70){$k_1$} 
%\put(110,17){$k_2$} 
% 

%
%\put(33,70){$P1$}

%
% 
% 
%\put(5,-12){ $\begin{array}{cc}
% \; & \begin{array}{ccc}
% L \hspace{1,8cm} &M\;\;\;\;\;\;\;\;\;&\;\;\;\;R\;\;\;\;\\
% \end{array} \\
% \begin{array}{c}
% T \\
% B \\
% \end{array} &
% \left( \begin{array}{ccc}
%  (-1,-1)\circlearrowleft & (-12,-11)^*& (-4,-7)^* \\
%  (-22,-12)^*& (3,-2)^*& (-9,-4)^* \\
%\end{array}\right)\\
%\end{array}  $}

%
%\put(-55,15){$(1,1)^*$}

%\put(-5,50){$T$}
%\put(70,50){$B$}

% 
%\put(35,65){\circle{2}}
% 
% 

% 
%\end{picture}
%\end{center}

\vspace{1cm}
\noindent There are 2 non absorbing states: $k_1$ (the initial state) and $k_2$, and 7 absorbing states. Player 1 has 2 actions: $T$ and $B$, and player 2 has 3 actions: $L$, $M$ and $R$. Playing $(B,L)$ in state $k_1$(or  $(T,L)$ in state $k_2$)  leads to state $k_2$.\\

%State $k_1$: $\begin{array}{cc}
% \; & \begin{array}{ccc}
% L \hspace{1,8cm} &M\;\;\;\;\;\;\;\;\;\;\;&\;\;\;\;R\;\;\;\;\\
% \end{array} \\
% \begin{array}{c}
% T \\
% B \\
% \end{array} &
% \left( \begin{array}{ccc}
%  (-1,-1)^{k_2}& (-10,-10)^*& (-10,-10)^* \\
%  (1,1)^*& (-10,-10)^*& (-10,-10)^* \\
%\end{array}\right)\\
%\end{array} $

%\vspace{0,5cm}

%State $k_2$: $\begin{array}{cc}
% \; & \begin{array}{ccc}
% L \hspace{1,8cm} &M\;\;\;\;\;\;\;\;\;&\;\;\;\;R\;\;\;\;\\
% \end{array} \\
% \begin{array}{c}
% T \\
% B \\
% \end{array} &
% \left( \begin{array}{ccc}
%  (-1,-1)\circlearrowleft & (-12,-11)^*& (-4,-7)^* \\
%  (-22,-12)^*& (3,-2)^*& (-9,-4)^* \\
%\end{array}\right)\\
%\end{array} $

%\vspace{0,5cm}

Consider the $\delta$-discounted game.    In any Nash equilibrium, Player 2 plays $L$ at stage 1.  There   exists a Nash equilibrium where $(T,L)$ is played at stage 1, and if  $\delta\geq 1/2$ there also exists a  Nash equilibrium where  $(B,L)$ is played at stage 1 and $(B,M)$ is played at stage 2. So for $\delta\geq 1/2$, $E_\delta$ contains the payoff $(1,1)$ as well as a payoff with  second  coordinate not greater than  -1 for player 2. The question is  now whether $E_\delta$ contains payoffs with second coordinate in $(-1,1)$, i.e. if at equilibrium player 1 can mix between $T$ and $B$ at stage 1.

\begin{lem} \label{lem2,5} In the game  {with initial state $k_2$} and discount $\delta$, the set of Nash (or perfect) equilibrium payoffs is:
$$\{(1-\delta^N)(-1,-1)+ \delta^N (u_1,u_2), N \in \N\cup\{+\infty\}, (u_1,u_2)\in \{(3,-2), (-6,-5)\} \}.$$
\end{lem}

\noindent{\bf Proof:} Let $(\sigma, \tau)$ be a Nash equilibrium of the stochastic game with initial state $k_2$. We prove that under $(\sigma, \tau)$,  either  $(T,L)$ is played at every stage, or $(T,L)$ is first played a certain  number $N$ of stages, then   followed if $N$ is finite  by an absorbing action profile with payoff $(3,-2)$ or $(-6,-5)$.

Denote by $\alpha'$ the best reply payoff of player 1 against the continuation strategy induced by $\tau$ after $(T,L)$ was played at stage 1. We have $\alpha'\in [-10,3]$, and define $\alpha=(1-\delta)(-1)+ \delta \alpha'\in [-10,3]$. Similarly, denote by $\beta'$ the best reply payoff of player 2 against the continuation strategy induced by $\sigma$ after $(T,L)$ was played at stage 1. And define   $\beta =(1-\delta)(-1)+ \delta \beta'\in [-13/2,-1]$. The strategies induced by $(\sigma, \tau)$ at stage 1 form a Nash equilibrium of the bimatrix game:

\vspace{0.5cm}
\begin{center}
$\begin{array}{cc}
 \; & \begin{array}{ccc}
 L \hspace{1,7cm} &M\;\;\;\;\;\;\;\;\;&\;\;\;\;R\;\;\;\;\\
 \end{array} \\
 \begin{array}{c}
 T \\
 B \\
 \end{array} &
 \left( \begin{array}{ccc}
  (\alpha, \beta) & (-12,-11)& (-4,-7) \\
  (-22,-12)& (3,-2)& (-9,-4) \\
\end{array}\right)\\
\end{array} $.
\end{center}
\vspace{0.5cm}

Let $x$, resp. $y_1$, $y_2$, $y_3$,  be the probability that $\sigma$, resp. $\tau$,  plays $T$, resp. $L$, $M$, $R$ at stage 1.
If $y_1y_2>0$, player 2 is indifferent between $L$ and $M$, which implies that $x=\frac{10}{21+\beta}\in [1/2,5/7]$. But then $R$ is strictly better than $M$ for player 2, hence a contradiction with $y_2>0$. If $y_1y_3>0$, then $y_2=0$ and $T$ is strictly better than $B$ for player 1, so $x=1$ which contradicts $y_3>0$. Consequently, if $y_1>0$   the only case is $y_1=1$, and $x=1$. $(\sigma, \tau)$ plays $(T,L)$ at stage 1. 

If $y_1=0$, we have a Nash equilibrium of the  game 

\vspace{0.5cm}
\begin{center}
$\begin{array}{cc}
 \; & \begin{array}{cc}
  M\;\;\;\;\;\;\;\;\;&\;\;\;\;R\;\;\;\;\\
 \end{array} \\
 \begin{array}{c}
 T \\
 B \\
 \end{array} &
 \left( \begin{array}{cc}
  (-12,-11)& (-4,-7) \\
  (3,-2)& (-9,-4) \\
\end{array}\right)\\
\end{array}$, 
\end{center}
\vspace{0.5cm}

and obtain\footnote{The profile $(T,R)$ is ruled  out by  action $L$ of player 2.} that either $(B,M)$, or the mixed action profile
\\
 $(1/3 T+ 2/3 B, 1/4 M+3/4 R)$, is played by $(\sigma, \tau)$ at stage 1. Notice that both $(B,M)$ and $(1/3 T+ 2/3 B, 1/4 M+3/4 R)$ lead with probability one to absorbing  states of  the stochastic game. 

Iterating  the argument from stage 1 on, leads to the inclusion of the Nash equilibrium payoffs set into  $\{(1-\delta^N)(-1,-1)+ \delta^N (u_1,u_2), N \in \N\cup\{+\infty\}, (u_1,u_2)\in \{(3,-2), (-6,-5)\} \}.$ It is then easy to see that  the conclusions of  Lemma \ref{lem2,5} hold.  

\vspace{1cm}

We now uniquely consider the stochastic game with initial state $k_1$. Define   $\Delta_1=\{\delta\in [0,1), \exists M\geq 1, (1-\delta^M)(-1)+3 \delta^M=1\}$, which is as in Example \ref{exa1}  the  countable  set of discount factors of the form $ {(\frac{1}{2})}^{1/M}$, with $M$ a positive integer.
 
Consider $\delta={(\frac{1}{2})}^{1/M} \in \Delta_1$, there exists a Nash equilibrium where $(B,L)$  is played at stage 1 having payoff $(1,u_2)$, with $u_2\leq -1$. Considering equilibria where player 1 mixes at stage 1, we obtain for $\delta$ in $\Delta_1$:
$$\{1\}\times [-1,1] \subset E'_\delta.$$

On the contrary, for $\delta\notin \Delta_1$, no Nash equilibrium of the stochastic game can mix between $T$ and $B$ at stage 1, so the equilibrium payoffs set satisfy:
$$E_\delta\cap \{(u_1,u_2), u_2\in (-1,1)\}=\emptyset.$$
 
\noindent This is enough to conclude that $E_\delta$ and $E'_\delta$ do not converge as $\delta$ goes to one. And one can easily check  that all arguments are robust to small perturbations of the payoffs of the stochastic game. \hfill $\Box$
  \end{exa}

\vspace{0,5cm}

We want to point out that the above examples are limited in several  ways. In particular:

 1) Many Folk theorems in the literature require the limit set to have non-empty interior. In both examples, $E_\delta=E'_\delta$ has empty interior for each discount factor.
 
 2) The examples are not robust to the introduction of an extensive-form    correlation device.  In example \ref{exa1}, if whenever $k_3$ is reached  the players can publicly observe the outcome of a fair coin tossing, they can correlate and play there   $(T,L)$ and $(B,R)$ with probability 1/2. With such correlation device, it is possible to obtain $(1/2,1/2)$ as an equilibrium payoff for all discount factors.  
 
   In example \ref{exa1.5},  the distance between the segment  $\{1\}\times [-1,1]$ and the set of normal-form  correlated equilibrium payoffs of the $\delta$-discounted game goes to 0 when $\delta$ goes to 1: consider a public correlation device mixing between the equilibrium with payoff $(1,1)$ and the equilibrium with payoff $(1-\delta^N)(-1,-1)+ \delta^N (3,-2)$ with $N$ the smallest integer such that $   \delta^N \geq 1/2$.
 
 3)  If the sets of equilibrium payoffs do not converge, some long-term equilibrium payoffs clearly emerge. 
  In example \ref{exa1}, for all  converging selections $(x_\delta)_\delta$ of $(E'_\delta)_\delta$, the limit payoff is $(1/2,0)$,  and for all discount $\delta$ the unique  stationary equilibrium payoffs is also $(1/2,0)$. Moreover, one can show that $E_\infty=\{(1/2,0)\}$, i.e. the unique uniform equilibrium payoff is $(1/2,0)$.   And playing $B$ in $k_1$ is somehow a risky option for player 1, since he can  immediately secure 1/2 by playing $T$ and has almost no chance to get a better payoff by playing $B$.  So even if the sets of equilibrium payoffs do not converge, the payoff $(1/2,0)$ clearly emerges, and we believe it can  be considered as the reasonable limit outcome of the stochastic game. If an expert is asked ``The game being played by  extremely patient players, which outcome is likely to form ?", we  would recommend the answer to be $(1/2,0)$. 
 
 In example \ref{exa1.5}, the limits of converging selections $(x_\delta)_\delta$ of $(E'_\delta)_\delta$ are the elements of the set 
 $E_\infty$  of uniform equilibrium payoffs of the stochastic game, which is the union of the singleton  $\{(1,1)\}$ and of the line segment joining $(1,-3/2)$ to $(3,-2)$. The limit set of stationary equilibrium payoffs $E$ defined by proposition \ref{pro1} is the pair $\{(1,1),(3,-2)\}$. $(1,1)$ is achieved by an equilibrium where player 2 plays $L$ in every period, and $(3,-2)$ is achieved by an equilibrium where player 1  first  plays $T$, then $B$ in each period. When the game is played by very patient players, we believe that one of these two payoff vectors is likely to occur.  \\
  
  The counterexample of the next section will not have these limiting properties and will be very robust in many aspects.  \\

 \section{Hidden Stochastic Games} \label{secHSG}
 
 We enlarge the model of stochastic games by assuming that at the beginning of every period, the players observe a public signal on the current state. We still  denote by  $K$, $I$ and $J$  respectively  the  finite sets of states,  actions for player 1 and actions for player 2, and we   introduce a finite set $S$ of public signals. As in the previous section,  $u_1$ and $u_2$ are the state  dependent utility functions from $K\times I \times J\longrightarrow \R$,  but now the transition function $q$ goes from $K \times I \times J$ to $\Delta(K\times S)$,   the set of probabilities over  $K\times S$, and there is an initial distribution $\pi$ in $\Delta(K\times S)$. The elements $K$, $I$, $J$, $S$, $u_1$, $u_2$, $q$ and $\pi$ are known to the players. 
 
 At   the first period,   a couple $(k_1,s_1)$ is selected according to $\pi$, and the players publicly observe $s_1$, but not $k_1$. The  players   simultaneously select actions $i_1\in I$ and $j_1\in J$, then  these actions are  publicly observed,   the stage payoffs are $u_1(k_1, i_1,j_1)$ for player 1 and $u_2(k_1, i_1, j_1)$ for player 2, and the play goes to period 2. At every period $t\geq 2$, a couple $(k_t,s_t)$ is selected according to  $q(k_{t-1},i_{t-1},j_{t-1})$, $k_t$ is the state of period $t$ but the players only observe the public signal $s_t$. Then they simultaneously select actions $i_t\in I$ and $j_t\in J$. These actions are  publicly observed, the stage payoffs are $u_1(k_t, i_t,j_t)$ for player 1 and $u_2(k_t, i_t,j_t)$ for player 2, and the play goes to the  period $t+1$.   
 Given a discount factor $\delta$ in $[0,1)$, the $\delta$-discounted hidden stochastic game   is the     game with   payoff functions    $(1-\delta) \sum_{t=1}^\infty \delta^{t-1} u_1(k_t,i_t,j_t)$ and $(1-\delta) \sum_{t=1}^\infty \delta^{t-1} u_2(k_t,i_t,j_t)$. We respectively  denote by  $E_\delta$ and $E'_\delta$  the sets of Nash equilibrium payoffs and sequential equilibrium payoffs of this game. 
 
 This is a generalization of the   model of stochastic game presented in section \ref{secSG}, where one has $S=K$ and $s_t=k_t$ for all $t$.  In the model of hidden stochastic game (HSG, for short), the players have   incomplete information on the current state, but this information is common to both players, and can be represented by a   belief $p_t$ on the state $k_t$. Given the initial signal $s_1$, the initial belief $p_1$  is the conditional probability induced by $\pi$ on $K$ given $s_1$.  The  belief $p_t$ is a random variable which can    be  computed\footnote{Notice that this belief does not depend on the {\it strategy} of the   players, as in repeated games with incomplete information,  but only on past actions played and public signals observed.} recursively from $p_{t-1}$ by Bayes' rule after observing the public signal $s_t$ and the  past actions $i_{t-1}$ and $j_{t-1}$.  We can thus associate to our HSG, an equivalent stochastic game  where the  state variable $p$ lies in $\Delta(K)$   and represents the common belief on the current state in the HSG, and where now actions and state variables are publicly observed, in addition to the public\footnote{In the equivalent stochastic game, the public signal $s$ gives no extra  information on past actions or on the state variable. Its unique  influence is that   it may be used by the players as a correlation device. Notice that the equivalent stochastic game is not a standard stochastic game as described in section \ref{secSG}.} signal $s$.  A strategy in the HSG uniquely defines an equivalent strategy in the stochastic game, and vice-versa. And in particular  the sets of equilibrium payoffs of the two games coincide. By definition, a stationary strategy in the associated stochastic game plays after every history a mixed action which only depends on the current state variable in $\Delta(K)$. And we will say that a strategy $\sigma$  in the HSG is stationary if the associated strategy in the stochastic game is stationary, that is if $\sigma$ plays after every history a mixed action which only depends on the current belief  in $\Delta(K)$.

Standard fixed-point (contraction) arguments show that $E_\delta$ and $E'_\delta$ are non-empty, and there exists a stationary equilibrium in the $\delta$-discounted associated stochastic game. We will also briefly consider  the set of uniform equilibrium payoffs $E_{\infty}$, defined as in the previous section.  

When there is a single player (for instance, when player 2 has a unique action), a hidden stochastic game is simply  a partially observable Markov decision process (POMDP), and if moreover player 1 plays constantly the same mixed action, we obtain a Hidden Markov model, which can be considered as the simplest model of dynamic Bayesian network. Hidden stochastic games generalize both standard stochastic games and POMDP. An interesting subclass of hidden stochastic games  is the following class of HSG  with known payoffs, where the public signals are rich enough for the players to know after every history what is the current payoff function. We write $q(k,i,j)(k',s)$ the probability in $[0,1]$ that the couple $(k',s)$ is chosen when the probability $q(k,i,j)$ is used.  

\begin{defi} \label{def1}  \rm The hidden stochastic game has {\it  known payoffs} if the set of states $K$ can be partitioned\footnote{We write $k_1\sim k_2$ whenever $k_1$ and $k_2$ are in the same equivalence class, or cell, of the partition.} in a way such that for all states $k$, $k'$, $k_1$, $k_2$, actions $i$, $i'$ in $I$, $j$, $j'$ in $J$, and signal $s$ in $S$:

1)if $k_1\sim k_2$ then $u_1(k_1,i,j)=u_1(k_2,i,j)$ and $u_2(k_1,i,j)=u_2(k_2,i,j)$  (two states in the same element of the partition induce the same payoff function), and

2)   if $q(k,i,j)(k_1,s)>0$ and $q(k',i',j')(k_2,s)>0$ then $k_1\sim k_2$ (observing the public signal is enough to deduce the element of the partition containing the current state).
\end{defi}

In a hidden stochastic game with known payoffs, the players know after every history the cell of the partition containing the current state, so when players choose their actions they know the current payoff function, as it happens in   a standard stochastic game . However they may not exactly    know the current state in $K$, so they are uncertain about the transition probabilities to the next state, and to the cell containing this state.  In a standard stochastic game, one can define: $k\sim k'$ if and only if $k=k'$, and the conditions of definition \ref{def1} are satisfied. Hence HSG with known payoffs generalize stochastic games, and this generalization is meaningful in several cases.

\begin{exa} \label{exa2} \rm The players are firms competing a la Cournot  on a market for a natural exhaustible resource. Only   two firms are present on this market, and in each period, each firm decides how much resource to extract (to produce). Then a price is set in order to equalize offer and demand, and all the production is sold at this price. Action sets are $I=\{0,...,M_1\}$ and $J=\{0,...,M_2\}$ where $M_f$ is the maximal possible production (e.g., in tonnes) of firm $f$.  The state variable $k$ is the amount of natural resources remaining (the stock),   and  the firms have incomplete information on remaining  stocks. They    have a common belief on the initial stock value $k_1$, and    there is a cap $M$ such that in each period,  if the current  stock $k$ is greater  than $M$ the firms just know that there are at least $M$ remaining resources, whereas if $k$ is at most $M$ the firms precisely know $k$.  Transitions are deterministic :  if in some state $k$,    actions $i$ and $j$ such that $i+j\leq k$ are played, then it is possible for the firms to actually produce the quantities $x$ and $y$, and the next state is\footnote{A more general variant for partially renewable  resources would read: the next state is$(1+r) (k-(i+j))$, where $r$ is the renewal rate.}   $k-(i+j)$. If $k> i+j$, the next state is 0 and the game is essentially over.  Payoffs are function of the actions, and  possibly  of the current state as well (when the state is lower than the sum of productions, or when the state does not exceed $M$ and the demand anticipates the scarcity of the resource).\end{exa}

\begin{exa} \label{exa3} \rm Consider an oligopoly with two firms on a market for a single good. In each period (e.g., a year) a firm chooses its selling price, as well as development and advertising budgets. The state variable $k$ represents the state of the market, which includes, but is not limited to, the current demand function for each firm, which is a function of the  current price profile. The state  also contains additional  information about  fundamentals   which will influence  the future evolution of the demand, such as the trend of the market, the development of close goods by other firms or the overall state of the economy.  In each period revenues are determined by the current demand function and the current prices chosen, and stage payoffs are the revenues minus   development and advertising budgets. Transitions of the state variable depend  on the state variable and the actions chosen, and  firms are   able to observe at the beginning of every period, at least  the current demand function but possibly  not all  characteristics of the state.
\end{exa}

Regarding limit equilibrium payoffs in hidden stochastic games, we know by proposition \ref{pro2} that there is no hope to obtain convergence of the sequences $(E_\delta)_\delta$ or $(E'_\delta)_\delta$. The following result shows that the situation is even more dramatic in our context of hidden stochastic games.

\begin{thm} \label{thm1} For each  $\varepsilon$ in $(0,\frac{5}{12}]$ and $r$ in $(0, \varepsilon/5)$, there exists a 2-player Hidden Stochastic Game   $\Gamma$ having the following properties:  
 \begin{enumerate}
\item
There are 13 states and  public signals,   four actions for each player, and  all payoffs lie in $[0,1]$,  
 
 \item  The game is symmetric between the players,  and  has known payoffs,  

\item
For all initial distributions and  discount factors,  the corresponding  set  of sequential equilibrium payoffs         contains a square of side $2r$, hence has full dimension,

\item There is an initial state $k_1$, perfectly known to the players, and there exist 
   two subsets  $\Delta_1$ and $\Delta_2$ of $[0,1)$, both containing discount factors arbitrarily close to 1,  such that: 

\begin{center}
for all $\delta$ in $\Delta_1$,   the set of sequential equilibrium payoffs $E'_{\delta}$ is the square $E_1$ centered in $(\epsilon,\epsilon)$ with side $2 r$, whereas 
for all $\delta$ in $\Delta_2$,   the set of sequential equilibrium payoffs  $E'_{\delta}$   is the square $E_2$ centered in $(1-\epsilon,1-\epsilon)$ with side $2 r$. 
\end{center}

Moreover  for $\delta$ in $\Delta_1\cup \Delta_2$,  the associated square  is also  the set of Nash equilibrium payoffs, the set of correlated equilibrium payoffs, and the set of communication equilibrium payoffs  of the $\delta$-discounted game, as well as the set of stationary equilibrium payoffs of the associated stochastic game with state variable the belief on the   states of the original game. 

There is no converging selection   of $(E_\delta)_\delta$, and $\Gamma$ has no uniform equilibrium payoff. 

\item The above conclusions   are robust to   perturbations of  the payoffs. Consider, for $\eta\in [0,\frac{r(\varepsilon-5r)}{2})$, a perturbed game $\Gamma(\eta)$ obtained by perturbing each  payoff of $\Gamma$ by at most $\eta$. The initial state being $k_1$,  denote by $E_\delta(\eta)$   the corresponding set of $\delta$-discounted Nash equilibrium payoffs. We have:
 \begin{eqnarray*}
\forall \delta \in  \Delta_1,&  
  E_\delta(\eta)\subset [\varepsilon-r- \eta, \varepsilon +r+  \eta]^2, \\   
\forall \delta \in  \Delta_2,&    E_\delta(\eta)\subset [1-\varepsilon-r- \eta, 1-\varepsilon +r+  \eta]^2.
\end{eqnarray*}
There is no converging selection     of $(E_\delta(\eta))_\delta$, and  $\Gamma(\eta)$ has no uniform equilibrium payoff.  
Finally, 
\begin{eqnarray*}
\lim_{\eta \to 0} \; \; \lim_{\delta \to 1, \delta \in \Delta_1}E'_\delta(\eta) =E_1& \; \rm{ and} \; &\lim_{\eta \to 0} \; \; \lim_{\delta \to 1, \delta \in \Delta_2}E'_\delta(\eta) =E_2,\\
 \lim_{\delta \to 1, \delta \in \Delta_1}  \limsup_{\eta \to 0} \; d(E_\delta(\eta),E_1)=0& \; \rm{ and} \; & \lim_{\delta \to 1, \delta \in \Delta_2}  \limsup_{\eta \to 0} \; d(E_\delta(\eta),E_2)=0.
 \end{eqnarray*}
\end{enumerate}\end{thm}
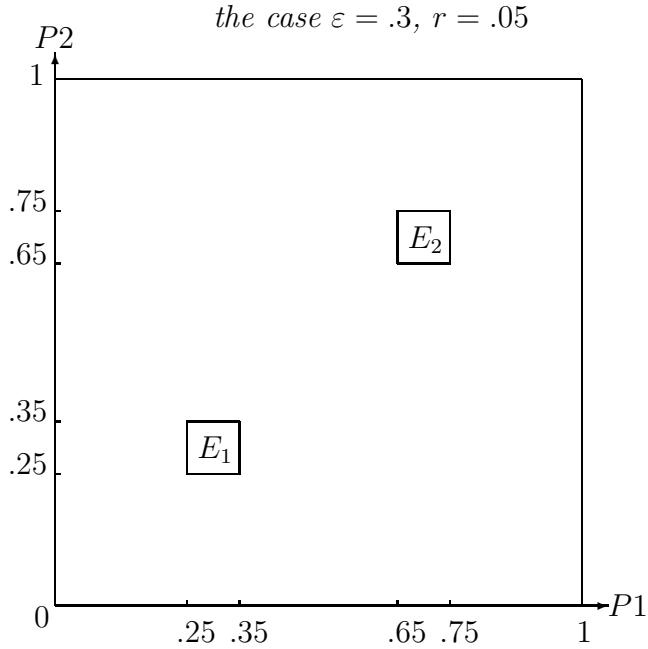
\unitlength 0,7mm
% sp�cifie l'unit�
%  \; \;{\rm and} \;\;    \lim_{\eta \to 0} E_\delta(\eta)=[\varepsilon-r, \varepsilon +r]^2, $
%   \;\; {\rm and} \;\;  \lim_{\eta \to 0} E_\delta(\eta)=[1-\varepsilon-r, 1-\varepsilon +r]^2.$

 \begin{center}

\begin{picture}(120,120)

 \put(30,110){\it the case $\varepsilon=.3$, $r=.05$}

 \put(0,0){\vector(1,0){105}}
 \put(0,0){\vector(0,1){105}}
  \put(100,0){\line(0,1){100}}
  \put(0,100){\line(1,0){100}}
  
   \put(-4,-4){0}
      \put(99,-7){1}
         \put(33,-7){$.35$}
            \put(23,-7){$.25$}
               \put(73,-7){$.75$}
                  \put(63,-7){$.65$}
          \put(35,0){\line(0,1){1}}
          \put(25,0){\line(0,1){1}}
          \put(75,0){\line(0,1){1}}
          \put(65,0){\line(0,1){1}}
          
           \put(-9,35){$.35$}

            \put(-9,25){$.25$}
               \put(-9,75){$.75$}
                  \put(-9,65){$.65$}
  \put(0,35){\line(1,0){1}}
          \put(0,25){\line(1,0){1}}
          \put(0,75){\line(1,0){1}}
          \put(0,65){\line(1,0){1}}

            \put(105,-2){$P1$}
             \put(-4,106){$P2$}
                   \put(-5,99){1}
                   
        \put(25,25){\line(0,1){10}}        
           \put(25,35){\line(1,0){10}}     
              \put(35,25){\line(0,1){10}}   
                 \put(25,25){\line(1,0){10}}     
                 
                         \put(65,65){\line(0,1){10}}        
           \put(65,75){\line(1,0){10}}     
              \put(75,75){\line(0,-1){10}}   
                 \put(65,65){\line(1,0){10}}
                 
            \put(27,28){$E_1$}       
              \put(67,68){$E_2$}  
\end{picture}
\end{center}
 
 \vspace{0,5cm}

The rest of the paper is devoted to the construction of the example  of theorem \ref{thm1}.  We progressively  introduce more and more ingredients in the construction, starting with a Markov chain on $[0,1]$, then a Markov Decision Process, then a zero-sum stochastic game with infinite state space,  a zero-sum HSG and finally our example\footnote{Apart the presentation, the  main differences with the 2013 example of Ziliotto are the following. Due to the zero-sum aspect in the  2013 example the game was non symmetric and equilibrium payoff sets had empty interior, this is taken care in the last construction of section \ref{secProof}.  In the non zero-sum case here the associated stochastic game with state variable in $\Delta(K)$ is not necessarily Nash-payoff equivalent to the HSG since public signals may be used in the HSG as correlation devices.  A main difference is that we need here equilibrium payoffs to go not only from 1/2 to 5/9, but from arbitrarily close to 0 to arbitrarily close to 1, so we    improve the construction by studying Markov chains and MDP with general parameters $\alpha$ and $\beta$ (which were equal  to 1/2 for both players  in the 2013 example). The asymmetry between players was obtained in 2013 by introducing a different structure for the MDP of player 2, whereas here the consideration of different parameters allows to stick to a symmetric, hence somehow simpler, construction.  We also   consider  non zero-sum perturbations of the payoffs, and have to deal with multiplicity of equilibria. Finally we also consider   multiple solution concepts: Nash, sequential, correlated and communication equilibria. }.

\section{Proof of Theorem \ref{thm1}}\label{secProof}

\subsection{A Markov chain on [0,1]} \label{sub1}

Given a parameter $\alpha\in (0,1)$, we consider the following Markov chain with state variable $q$ in $[0,1]$ and  initial state $q_0=1$.  Time is discrete, and if $q_t$ is the state of period $t$ then with probability $\alpha$ the next state $q_{t+1}$ is $\alpha \,q_t$ and with probability $1-\alpha$ the next state $q_{t+1}$ is 1.

 
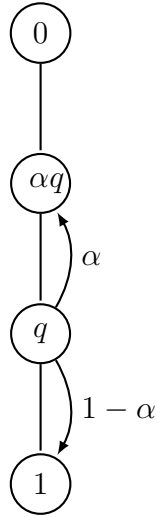
\begin{center}

\begin{tikzpicture}[>=stealth',shorten >=1pt,auto,node distance=4cm,thick,main
 node/.style={circle,draw,font=\Large\bfseries}]

\node [draw,text width=0.3cm,text centered,circle] (O) at (0,0) {$0$};
\node [draw,text width=0.3cm,text centered,circle] (A) at (0,-2) {$\alpha q$};
\node [draw,text width=0.3cm,text centered,circle] (B) at (0,-4) {$q$};
\node [draw,text width=0.3cm,text centered,circle] (C) at (0,-6) {$1$};

\draw[->,>=latex] (B) to[bend right] node[midway, right] {$\alpha$}(A);
\draw[->,>=latex] (B) to[bend left] node[midway,right] {$1-\alpha$}(C);
\draw[-,>=latex] (O) edge[] node[midway,above right] {}(A);
\draw[-,>=latex] (A) edge[] node[midway,above right] {}(B);
\draw[-,>=latex] (B) edge[] node[midway,above right] {}(C);

\end{tikzpicture}
\end{center}
Because of the transitions, the set of states that can be reached   is the countable set $\{\alpha^a, a\in \N\}$. This Markov chain can be viewed as follows: there is an infinite sequence $X_1$, ..., $X_t$, ...  of i.i.d. Bernouilli random variables with success parameter $\alpha$, we add an initial constant variable $X_0=0$, and at any  period $t$ the  state of the Markov chain is $\alpha^a$ if and only if the last $a$ (but not $a+1$) realizations of the Bernouilli variables have been successful, i.e. iff $X_{t-a}=0$ and $X_{t'}=1$  for $t-a+1\leq t'\leq t$. 

In the next  subsection, the variable $q$ will be interpreted as a {\it risk} variable with the following interpretation. Suppose a decision-maker observes the realizations of the Markov chain, and has to decide as a function of $q$ when he will take a risky action, having  probability  of success $1-q$ and probability of failure $q$. He would like $q$ to be as small as possible, but time is costly and there is a discount factor $\delta$. For $a$ in $\N$, we denote by $T_a$ the stopping time of the first period where the risk is $\alpha^a$, i.e.
$$T_a=\inf\{t\geq 1, q_t\leq \alpha^a\}.$$
\noindent If $a=0$, then $T_a=1$ and $\delta^{T_a}=\delta$. If $a\geq 1$, then $T_a$ is a random variable which law can be easily computed by induction. Indeed, we have: 
\begin{eqnarray*} T_a &= &  {T_{a-1}}+1+ \1_{X_{T_a}=0}\; T'_a,
 \end{eqnarray*}
\noindent where $T'_a$ has the same law as $T_a$ and is independent from $X_{T_a-1}$. As a consequence, $$\E(T_a)= \frac{1}{\alpha} (1+\E(T_{a-1})).$$
\noindent   $\E(T_a)$ grows exponentially with $a$, and this is  an important feature of our  counterexample: while slightly decreasing the risk $\alpha^a$ in the bounded set $(0,1]$, the number of stages one may have  to wait before reaching the new risk level greatly increases. 

The expectation of $\delta^{T_a}$ will play an important role in the sequel and can be easily computed as well (see e.g. lemma 2.2 and proposition 2.6   in \cite{PG83}).
 \begin{lem} \label{lem1.5}
$$\E( \delta^{T_a})=\frac{1- \alpha \delta}{1-\alpha+(1-\delta) \alpha^{-a} \delta^{-a-1}}.$$
\end{lem}

\subsection{A Markov Decision Process on [0,1]} \label{sub2}

We introduce a player  who observes the realizations of the above Markov chain and can choose as a function of the state $q$ when he will take a risky action, having probability  of success $1-q$ and probability of failure $q$. In case of success, the payoff of the player will be $R$ at all subsequent stages, where $R$ is a fixed positive reward. The payoff is 0 at any stage before taking the risky action,   and at any stage after the risky action has been taken unsuccessfully.  Overall payoffs are discounted with discount $\delta$.

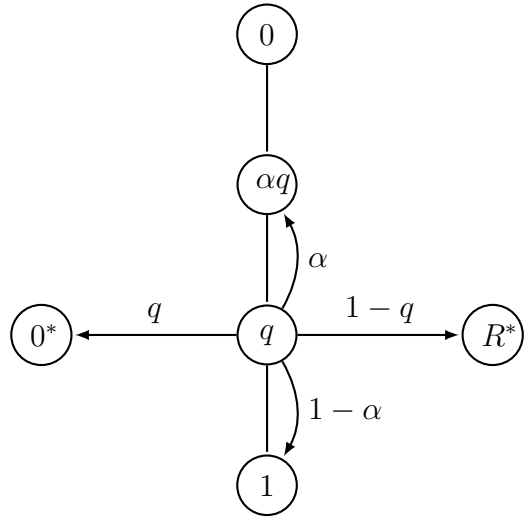
\begin{center}

\begin{tikzpicture}[>=stealth',shorten >=1pt,auto,node distance=4cm,thick,main
 node/.style={circle,draw,font=\Large\bfseries}]
% \draw (0,8) node[below]{\textbf{Player 1}};
% \draw (9,8) node[below]{\textbf{Player 2}};
\node [draw,text width=0.3cm,text centered,circle] (O) at (0,0) {$0$};
\node [draw,text width=0.3cm,text centered,circle] (A) at (0,-2) {$\alpha q$};
\node [draw,text width=0.3cm,text centered,circle] (B) at (0,-4) {$q$};
\node [draw,text width=0.3cm,text centered,circle] (C) at (0,-6) {$1$};
\node [draw,text width=0.3cm,text centered,circle] (D) at (-3,-4) {$0^*$};
\node [draw,text width=0.3cm,text centered,circle] (E) at (3,-4) {$R^*$};

\draw[->,>=latex] (B) to[bend right] node[midway, right] {$\alpha$}(A);
\draw[->,>=latex] (B) to[bend left] node[midway,right] {$1-\alpha$}(C);
\draw[-,>=latex] (O) edge[] node[midway,above right] {}(A);
\draw[-,>=latex] (A) edge[] node[midway,above right] {}(B);
\draw[-,>=latex] (B) edge[] node[midway,above right] {}(C);
\draw[->,>=latex] (B) edge[] node[midway, above] {$q$}(D);
\draw[->,>=latex] (B) to[] node[midway, above] {$1-q$} (E);

%\draw[->,>=latex] (O) to[] node[midway, right] {$J_1,s_0^*$}(B);

\end{tikzpicture}
\end{center}

In this MDP with finite actions set, there exists  a pure stationary optimal strategy. Notice that a pure stationary strategy of the player can be represented by a non negative integer $a$, corresponding  to the risk threshold $\alpha^a$. We define the $a$-strategy of the player as the strategy where he takes the risky action as soon as the state variable of the Markov chain does not exceed $\alpha^a$. The expected discounted payoff induced is  
$$\E \left( (1-\delta^{T_a}) 0 + \delta^{T_a}( \alpha^a 0+ (1-\alpha^a)R)\right)= R \, (1-\alpha^a) \, \E(\delta^{T_a}).$$
Hence using lemma \ref{lem1.5}, we obtain:

\begin{lem} \label{lem2.5}
The payoff of the $a$-strategy in the MDP with parameter $\alpha$ and discount $\delta$ is:
$$  \frac{(1-\alpha^a)(1- \alpha \delta)R}{1-\alpha+(1-\delta) \alpha^{-a} \delta^{-a-1}}.$$
\end{lem}

This payoff is proportional to $R>0$, hence the optimal strategies  do not depend on the value of $R$. Intuitively this is clear,  counting the reward in Dollars or Euros does not affect the strategic problem of the decision-maker. This problem is now to choose a non negative integer $a$ maximizing the above payoff function.

\begin{defi}\label{def2} Define, for all $a$ in $\R_+$, 
$$s_{\alpha, \delta}(a)=(1-\alpha^a) \E(\delta^{T_a})= \frac{(1-\alpha^a)(1- \alpha \delta)}{1-\alpha+(1-\delta) \alpha^{-a} \delta^{-a-1}},$$
\noindent and let   $v_{\alpha,\delta}=\max_{a \in \N} s_{\alpha, \delta}(a)$  denote the value of the $\delta$-discounted MDP with parameter $\alpha$ and reward $R=1$ .\end{defi}

\noindent $v_{\alpha,\delta}=\max_{a \in \N}   s_{\alpha, \delta}(a)$ is clear\footnote{One can verify analytically that the maximum of $s_{\alpha, \delta}$ over $\N$ is achieved, since  $0=s_{\alpha, \delta}(0)=\lim_{+\infty}s_{\alpha, \delta}$.} since there exists a pure optimal stationary strategy in the $\delta$-discounted MDP. The parameter $\alpha$ being fixed, we are now interested in maximizing $s_{\alpha, \delta}$ for $\delta$ close to 1.  Differentiating the   function $(a\mapsto \frac{(1-\alpha^a)}{1-\alpha+(1-\delta) \alpha^{-a}})$
and proceeding by asymptotical equivalence when $\delta$ goes to 1, naturally  leads to the introduction of the following quantity.

\begin{defi}\label{def3} When  $\delta \in [\alpha,1)$, we define $a^*= a^*(\alpha, \delta)$ in $\R_+$  such that: $$\alpha^{a^*}=\sqrt{\frac{1-\delta}{1-\alpha}}.$$
Let  $\Delta_1(\alpha)=\{1-(1-\alpha)\alpha^{2a}, a\in \N\}$   be the set of discount factors $\delta$ such that $a^*(\alpha, \delta)$ is an integer, and   let $\Delta_2(\alpha)=\{1-(1-\alpha)\alpha^{2a+\eta}, a\in \N, \eta \in [-3/2,3/2] \}$ be the set of discount factors $\delta$ such that $a^*(\alpha, \delta) \in \m{N}+[1/4,3/4]$. 
\end{defi}

\noindent $\Delta_1(\alpha)$ and $\Delta_2(\alpha)$ contain discount factors arbitrarily close to 1. $a^*$ can be expressed in closed form as  $a^*=\frac{\ln(1-\delta)-\ln(1-\alpha)}{2\ln \alpha}$. Since $\delta^{\ln(1-\delta)}$ converges to 1 when $\delta$ goes to 1,  we obtain: $\delta^{a^*}\xrightarrow[\delta \to 1] {}1.$  

\begin{pro}\label{pro3}  $\;$

1) $v_{\alpha,\delta}\xrightarrow[\delta \to 1] {}1.$  

2) For $\alpha<1/4$ and $\delta\in \Delta_1(\alpha)$, the $a^*(\alpha, \delta)$-strategy is optimal in the MDP and 
$$\lim_{\delta \to 1, \delta \in \Delta_1(\alpha)}\; \frac{1-v_{\alpha,\delta}}{ \sqrt{ {1-\delta} }}=\frac{2}{\sqrt{1-\alpha}}.$$

3) For all $\alpha$, $$\liminf_{\delta \to 1, \delta \in \Delta_2(\alpha)} \; \frac{1-v_{\alpha,\delta}}{  \sqrt{1-\delta}}\geq \frac{1}{\sqrt{\alpha^{1/2} (1-\alpha)}}.$$
\end{pro}
The convergence property in 1) is very intuitive: when $\delta$ is high, the decision-maker can wait for the state variable to be very low, so that she takes the risky action with high probability of success. Points 2) (when $\alpha<1/4$) and 3) give asymptotic expansions for the value $v_{\alpha,\delta}$ when $\delta$ goes to 1, respectively of the form 
$v_{\alpha,\delta}=1-2 \sqrt{\frac{1-\delta}{1-\alpha}}+ \sqrt{{1-\delta}}\,  \varepsilon_{\alpha}(\delta)$ and $v_{\alpha,\delta}\leq 1- \sqrt{\frac{1-\delta}{\alpha^{1/2}(1-\alpha)}}+ \sqrt{1-\delta}\,  \varepsilon'_{\alpha}(\delta)$, where $\varepsilon_{\alpha}$ and $\varepsilon'_{\alpha}$ are functions with limit 0 when $\delta$ goes to 1. Later on,  the parameter $\alpha$ will be small, and the situation of the associated  player  will be much better when $\delta$ is close to 1 in $\Delta_1(\alpha)$ compared to when $\delta$ is close to 1 in $\Delta_2(\alpha)$. The proof of  proposition \ref{pro3} is based on simple computations that are presented in the Appendix.

\subsection{A zero-sum stochastic game with perfect information} \label{sub3}

We fix here two parameters $\alpha$ and $\beta$ in $(0,1)$, and define a 2-player zero-sum stochastic game $\Gamma_{\alpha, \beta}$ with infinite state space: 
 $$X=\{(1,q), q \in [0,1]\}\;  \cup \;  \{(2,l), l \in [0,1]\} \; \cup \;  0^* \; \cup \; 1^*.$$

\begin{center}
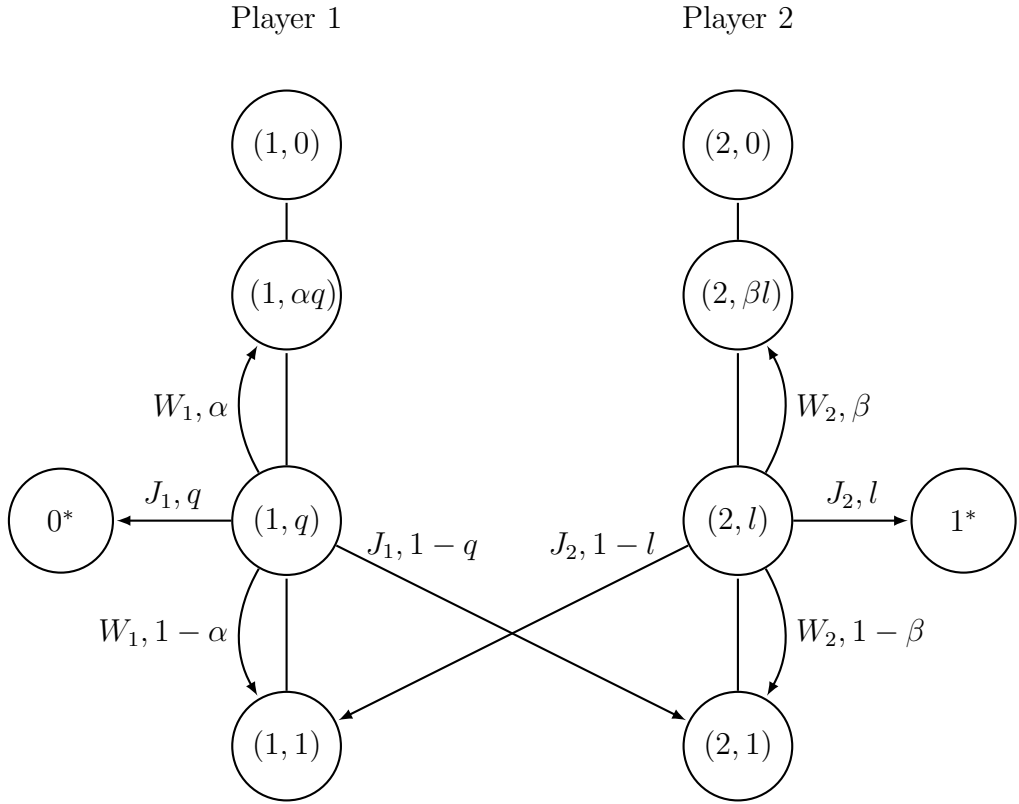
\begin{figure}[!h]

\begin{center}
\begin{tikzpicture}[>=stealth',shorten >=1pt,auto,node distance=4cm,thick,main
 node/.style={circle,draw,font=\Large\bfseries}]

\draw (0,10) node[below]{\text{Player 1}};
\draw (6,10) node[below]{\text{Player 2}};

\draw (1.8,3) node[below]{$J_1,1-q$};
\draw (4.2,3) node[below]{$J_2,1-l$};

\node [draw,text width=1cm, text centered,circle] (O) at (0,0) {$(1,1)$};
\node [draw,text width=1cm,text centered,circle] (A) at (0,3) {$(1,q)$};
\node [draw,text width=1cm,text centered,circle] (B) at (0,6) {$(1,\alpha q)$};
\node [draw,text width=1cm,text centered,circle] (C) at (0,8) {$(1,0)$};

\node [draw,text width=1cm,text centered,circle] (D) at (-3,3) {$0^*$};

\node [draw,text width=1cm,text centered,circle] (O') at (6,0) {$(2,1)$};
\node [draw,text width=1cm,text centered,circle] (A') at (6,3) {$(2,l)$};
\node [draw,text width=1cm,text centered,circle] (B') at (6,6) {$(2,\beta l)$};
\node [draw,text width=1cm,text centered,circle] (C') at (6,8) {$(2,0)$};

\node [draw,text width=1cm,text centered,circle] (D') at (9,3) {$1^*$};

\draw[->,>=latex] (A) to[] node[midway, above] {$J_1,q$}(D);

\draw[->,>=latex] (A) to[bend left] node[midway, left] {$W_1,\alpha$}(B);
\draw[->,>=latex] (A) to[bend right] node[midway, left] {$W_1,1-\alpha$}(O);

%\draw[->,>=latex] (A) to[] node[very near start, above right] {$J_1,1-q$}(O');
\draw[->,>=latex] (A) to[] node[]{}(O');

\draw[-,>=latex] (O) to[] node[] {}(A);
\draw[-,>=latex] (A) to[] node[] {}(B);
\draw[-,>=latex] (B) to[] node[] {}(C);

\draw[->,>=latex] (A') to[] node[midway, above] {$J_2,l$}(D');
\draw[->,>=latex] (A') to[bend right] node[midway, right] {$W_2,\beta$}(B');
\draw[->,>=latex] (A') to[bend left] node[midway, right] {$W_2,1-\beta$}(O');

%\draw[->,>=latex] (A') to[] node[very near start, above left] {$J_2,1-l$}(O);
\draw[->,>=latex] (A') to[] node[] {}(O);

\draw[-,>=latex] (O') to[] node[] {}(A');
\draw[-,>=latex] (A') to[] node[] {}(B');
\draw[-,>=latex] (B') to[] node[] {}(C');

\end{tikzpicture}
 
\vspace{0,5cm}

\caption{The stochastic game $\Gamma_{\alpha,\beta}$}
\end{center}
\end{figure}
\end{center}

\vspace{0,5cm}

\noindent    The initial state is $(2,1)$.  The sum of the payoffs of the players is constant\footnote{so strictly speaking, the game is constant-sum and not zero-sum, but we make
     the usual language abuse.} equal to 1. States $0^*$ and $1^*$ are absorbing states with, respectively, payoffs 0 and 1 to player 1. The payoffs only depend on the states, and the payoff of player 1 is 0 in a state of the form $(1,q)$ , and 1 in a state of the form $(2,l)$. Each player has 2 actions: Wait    or Jump. Transitions in a state $(1,q)$ are controlled by player 1 only: if player 1 Waits in state $(1,q)$, then the next state is $(1, \alpha q)$ with probability $\alpha$ and $(1,1)$ with probability $1-\alpha$, as in the MDP of subsection \ref{sub2}, and if player 1 Jumps in state $(1,q)$, then the next state is $0^*$ with probability $q$ and $(2,1)$ with probability $1-q$. Similarly,  transitions in a state $(2,l)$ are controlled by player 2 only: if player 2 Waits in state $(2,l)$, then the next state is $(2, \beta l)$ with probability $\beta$ and $(2,1)$ with probability $1-\beta$, and if player 2 Jumps in state $(2,l)$, then the next state is $1^*$ with probability $l$ and $(1,1)$ with probability $1-l$. Payoffs are discounted with discount factor $\delta\in [0,1)$, and the value of the stochastic game is denoted $v_{\alpha,\beta, \delta}$.
     
    The strategic aspects of this game have strong similarities with those of the previous MDP. Consider for instance Player 1, his payoff is 0 in $0^*$ and all states $(1,q)$, and his payoff is 1 in $1^*$ and the states $(2,l)$. Starting from  state (1,1), the only possibility for Player 1 to obtain positive payoffs is to Jump at some period to try to reach the state $(2,1)$. He can wait for the state to be $(1,q)$ with $q$ small, so that the  risk of reaching the state $0^*$ while jumping  is low, but each period in a state $(1,q)$ gives him a null payoff  so   he should not wait too long. The situation is symmetric for player 2, apart from the fact that the initial state is $(2,1)$, hence controlled by him. 
    
    Since the game is discounted and states are controlled by a single player, it is natural to look at pure stationary\footnote{Notice that $a$ and $b$-strategies are not fully defined in definition \ref{def4}, since they do not specify the actions played in the absorbing states nor in the states controlled by the other player. Since these actions have no impact on the game, we will simply ignore them.}  strategies of the players.
    
\begin{defi} \label{def4}    For $a$    in $\N$,
the  $a$-strategy of Player 1  is the strategy where Player 1  Jumps in a state $(1,q)$ if and only if  $q\leq \alpha^a$. Similarly, for $b$ in $\N$ the  $b$-strategy of Player 2  is the strategy where Player 2  Jumps in a state $(2,l)$ if and only if  $l\leq \beta^b$.  And we denote by  $g_{\alpha,\beta,\delta}(a,b)$ the payoff of Player 1 in the stochastic game where Player 1 uses the $a$-strategy and Player 2 uses the $b$-strategy. 
 \end{defi}
 
Assume   that Player 2 uses a $b$-strategy. Then Player 1 faces a MDP with finite action sets, hence he/she has a pure stationary best reply, that is Player 1 has a best reply in the stochastic game in the form of a $a$-strategy. Similarly, if Player 1 uses a $a$-strategy, Player 2 has a best reply in the stochastic game in the form of a $b$-strategy. It is then  natural to consider the game restricted to $a$- and $b$-strategies.

\begin{lem} \label{lem4} For  $a$ and $b$ in $\N$,
$$g_{\alpha,\beta,\delta}(a,b)=\frac{1-s_{\beta, \delta}(b)}{1- s_{\alpha, \delta}(a) s_{\beta, \delta}(b)}.$$
\end{lem}
     
\noindent{\bf Proof:} Recall that $s_{\alpha, \delta}(a)=(1-\alpha^a) \E_\alpha(\delta^{T_a})$, where $T_a$ is the random variable defined in subsection \ref{sub1} and $\E_\alpha$ denotes the expectation for the Markov chain with parameter $\alpha$. Similarly,  one has   $s_{\beta, \delta}(b)=(1-\beta^b) \E_\beta(\delta^{T_b})$.

Starting from the initial state, with probability $\beta^b$ the first Jump of player 2 will end up in $1^*$ and the payoff for player 1 will be 1 in each period, and with probability $1-\beta^b$ the game will first stay $T_b$ stages in a state controlled by player 2  and then reach the state $(1,1)$. This gives:
$$g_{\alpha,\beta,\delta}(a,b)= \beta^b+(1-\beta^b)  \E_\beta \left(  (1-\delta^{T_b}) +  \delta^{T_b} g'_{\alpha,\beta,\delta}(a,b)\right),$$
\noindent where $g'_{\alpha,\beta,\delta}(a,b)$ denotes the payoff of the $a$-strategy against the $b$-strategy in the game with initial state $(1,1)$. So $g_{\alpha,\beta,\delta}(a,b)=
 1+ s_{\beta, \delta}(b) (-1+g'_{\alpha,\beta,\delta}(a,b))$. Similarly, 
 $$g'_{,\alpha,\beta,\delta}(a,b)= \alpha^a 0 +(1-\alpha^a)  \E_\alpha (\delta^{T_a})   g_{\alpha,\beta,\delta}(a,b),$$
\noindent so $g'_{\alpha,\beta,\delta}(a,b)=s_{\alpha, \delta}(a) g_{\alpha,\beta,\delta}(a,b).$ Hence the result of lemma \ref{lem4}.\\

 Let us  come back to the consideration that Player 2 plays a $b$-strategy, and denote by $R$ the best payoff that Player 1 can obtain against  this strategy from the state $(1,1)$ (if the play  never reaches this state, then player 1 has nothing to do and gets  a payoff of 1 in each period). We have seen that Player 1 has a best reply in the form  of a $a$-strategy, and finding the best $a$ is equivalent to finding a pure optimal strategy in the MDP of subsection \ref{sub2} with reward $R$. But we have seen in  subsection \ref{sub2} that this optimal value for $a$ does not depend on $R$, and simply maximizes $s_{\alpha, \delta}(a)$. This implies that the best reply of player 1 does not depend on $b$, and the corresponding $a$-strategy is a {\it dominant} strategy of player 1 in the zero-sum stochastic game restricted to pure stationary strategies. The existence of dominant strategies in a zero-sum game is rather rare, and this is an important property of the present example. It can be verified   analytically  by looking at the function $g_{\alpha,\beta,\delta}$:  for all $b$, it is increasing in $s_{\alpha,\delta}(a)$, and for all $a$, it is decreasing in $s_{\beta,\delta}(b)$. This proves $1)$  in the  proposition below. 
 
 \begin{pro} \label{pro4} Let $a^\#$ and $b^\#$  be   respectively maximizers of   $s_{\alpha,\delta}(a)$ for $a$ in $\N$,  and of $s_{\beta,\delta}(b)$ for $b$ in $\N$, i.e. be non negative integers such that $s_{\alpha,\delta}(a^\#)=v_{\alpha, \delta}$ and $s_{\beta,\delta}(b^\#)=v_{\beta, \delta}$.
 
 1) The $a^\#$-strategy, resp. the $b^\#$-strategy,     is a dominant strategy for player 1, resp. player 2,  in the zero-sum stochastic game restricted to pure stationary strategies.
 
 2)  The $a^\#$-strategy, resp. the $b^\#$-strategy,      is an optimal  strategy for player 1, resp. player 2,  in the zero-sum stochastic game $\Gamma_{\alpha, \beta}$.
 
 3) The value of $\Gamma_{\alpha, \beta}$ satisfies:
 
$$v_{\alpha,\beta, \delta}=\frac{1-v_{\beta, \delta} }{1- v_{\alpha, \delta}  v_{\beta, \delta}}.$$
 \end{pro}
 
 \noindent{\bf Proof:}  $2)$ The strategy profile induced by $(a^\#,b^\#)$ is a Nash equilibrium of the game  $\Gamma_{\alpha, \beta}$  restricted to pure stationary strategies. Since   
  against a pure stationary strategy each player has a pure stationary best reply, this strategy profile is indeed a Nash equilibrium of the   game $\Gamma_{\alpha, \beta}$. Hence the value of $\Gamma_{\alpha, \beta}$ is the payoff induced by this strategy profile, and $3)$ follows. \\
  
Notice that $v_{\alpha, \alpha, \delta}=\frac{1}{1+v_{\alpha, \delta}}\xrightarrow[\delta \to 1]{} \frac{1}{2}$.  We are interested in cases where $\alpha\neq \beta$, and  the next proposition is  a building brick for our global construction.
  
   \begin{pro} \label{pro5} For each $\varepsilon>0$, there exists $n_0 \in \m{N}^*$ such that for all $n \geq n_0$, and $\alpha:=1/n$ and $\beta:=1/(n+1)$, we have:
   $$\limsup_{\delta\to 1}v_{\alpha,\beta, \delta}\geq 1-\varepsilon,\; {\rm and}\; \liminf_{\delta\to 1}v_{\alpha,\beta, \delta}\leq \varepsilon.$$
    \end{pro}
    
     \noindent{\bf Proof:} We proceed in $2$ steps.
     
       \underline{Step 1:} Define $\Delta_1(\alpha,\beta):=\Delta_1(\alpha)\cap \Delta_2(\beta)$, that is:
       \begin{equation*}
      \Delta_1(\alpha,\beta) =  \{\delta \in [0,1), \exists (a,b,\eta) \in \N^2 \times [-3/2,3/2],\delta=1-(1-\alpha)\alpha^{2a} =1-(1-\beta)\beta^{2b+\eta}\}.
      \end{equation*}
       \noindent Discount factors in $\Delta_1(\alpha,\beta)$ simultaneously favor  player 1 and disfavor  player 2 in their respective MDP: for $\delta\in \Delta_1(\alpha,\beta)$,  we have by proposition \ref{pro3} that $v_{\alpha,\delta}=1-2 \sqrt{\frac{1-\delta}{1-\alpha}}+ \sqrt{{1-\delta}}\,  \varepsilon_{\alpha}(\delta)$ and $v_{\beta,\delta}\leq 1- \sqrt{\frac{1-\delta}{\beta^{1/2}(1-\beta)}}+ \sqrt{1-\delta}\,  \varepsilon'_{\beta}(\delta)$, with  $\lim_{\delta \to 1} \varepsilon_{\alpha}=\lim_{\delta \to 1} \varepsilon'_{\beta}=0$.     Since $v_{\alpha,\beta, \delta}=\frac{1-v_{\beta, \delta} }{1- v_{\alpha, \delta}  v_{\beta, \delta}}$ is decreasing in $v_{\beta, \delta}$, we obtain:
     \begin{eqnarray*}
     v_{\alpha,\beta, \delta}&\geq &\frac{\sqrt{\frac{1-\delta}{\beta^{1/2}(1-\beta)}}- \sqrt{1-\delta}\,  \varepsilon'_{\beta}(\delta)}{1- \left(1-2 \sqrt{\frac{1-\delta}{1-\alpha}}+ \sqrt{{1-\delta}}\,  \varepsilon_{\alpha}(\delta)\right)\left( 1- \sqrt{\frac{1-\delta}{\beta^{1/2}(1-\beta)}}+ \sqrt{1-\delta}\,  \varepsilon'_{\beta}(\delta)\right)},\\
 &    \geq &\frac{\sqrt{\frac{1-\delta}{\beta^{1/2}(1-\beta)}}- \sqrt{1-\delta}\,  \varepsilon'_{\beta}(\delta)}
     { \sqrt{\frac{1-\delta}{\beta^{1/2}(1-\beta)}}+2 \sqrt{\frac{1-\delta}{1-\alpha}}+ \sqrt{1-\delta}\,  \varepsilon''(\delta)},\; {\rm where}\; \lim_{\delta \to 1} \varepsilon''=0.
     \end{eqnarray*}
     \noindent This implies,  if $\Delta_1(\alpha,\beta)$ contains discount factors arbitrarily close to 1:
      \begin{equation} \label{eq3}\liminf_{\delta \to 1, \delta \in \Delta_1(\alpha,\beta)}  v_{\alpha,\beta, \delta} \geq \frac{1}{1+2 \sqrt{\frac{\beta^{1/2}(1-\beta)}{1-\alpha}}}.\end{equation}
     In the same vein, we define $\Delta_2(\alpha,\beta):=\Delta_2(\alpha)\cap \Delta_1(\beta)$, that is:
     \begin{equation*}
      \Delta_2(\alpha,\beta)=\{\delta \in [0,1), \exists (a,b,\eta) \in \N^2 \times [-3/2,3/2], \delta=1-(1-\alpha)\alpha^{2a+\eta} =1-(1-\beta)\beta^{2b}\}.
      \end{equation*}
       \noindent Discount factors in $\Delta_2(\alpha,\beta)$ simultaneously  disfavor  player 1 and  favor  player 2 in their respective MDP, and similar computations as above show  that  if $\Delta_2(\alpha,\beta)$ contains discount factors arbitrarily close to 1,
      \begin{equation} \label{eq4} \limsup_{\delta \to 1, \delta \in \Delta_2(\alpha,\beta)}  v_{\alpha,\beta, \delta} \leq \frac{1}{1+\frac{1}{2} \sqrt{\frac{(1-\beta)}{\alpha^{1/2}(1-\alpha)}}}.\end{equation}
 Our goal, inspired by (\ref{eq3}) and (\ref{eq4}),  is now to prove that there exists $\alpha$ and $\beta$ arbitrarily small such that both  $\Delta_1(\alpha,\beta)$ and $\Delta_2(\alpha,\beta)$ contain discount factors arbitrarily close to 1. \\
 
\underline{Step 2:} 

   % This implies that both $\Delta_1(\alpha, \beta)$ and $\Delta_2(\alpha, \beta)$ contain discount factors arbitrarily close to 1. 
   We want to prove that for $n$ big enough, there exists an infinite number of pairs $(a,b,\eta) \in \m{N}^2\times [-3/2,3/2]$ verifying
   \begin{equation*}
   1-(1-\alpha) \alpha^{2a}=1-(1-\beta) \beta^{2b+\eta},
   \end{equation*}
   that is,
   \begin{equation*}
   \ln(\beta)^{-1}\left[\ln((1-\alpha)/(1-\beta))+2a \ln(\alpha) \right]=2b+\eta.
   \end{equation*}
   Let $A(\alpha,\beta):=\ln(\beta)^{-1} \ln((1-\alpha)/(1-\beta))$ and $B(\alpha,\beta):=\ln(\beta)^{-1} \ln(\alpha)-1$. 
   %For $A(\alpha,\beta)$ and $B(\alpha,\beta)$ small enough, the equation
   The last equation can be written as
   \begin{equation*}
   A(\alpha,\beta)+2B(\alpha,\beta) a=2(b-a)+\eta.
   \end{equation*}
   If $B(\alpha,\beta)<1/4$, then this equation has an infinite number of solutions $(a,b,\eta) \in \m{N}^2 \times [-3/2,3/2]$. Set $\alpha_n:=1/n$ and $\beta_n:=1/(n+1)$. For $n$ big enough, we have $B(\alpha_n,\beta_n)<1/4$. This implies that $\Delta_1(\alpha_n, \beta_n)$ contains discount factors arbitrarily close to 1, and the proof is similar for $\Delta_2(\alpha_n, \beta_n)$. 
   \\
Let $\epsilon$ and $n_0 \in \m{N}$ such that for all $n \geq n_0$, both $\Delta_1(\alpha_n, \beta_n)$ and $\Delta_2(\alpha_n, \beta_n)$ contain discount factors arbitrarily close to 1, and $\left(1+2 \sqrt{\frac{\beta_n^{1/2}(1-\beta_n)}{1-\alpha_n}}\right)^{-1} \geq 1-\epsilon$, and $\left(1+\frac{1}{2} \sqrt{\frac{(1-\beta_n)}{\alpha_n^{1/2}(1-\alpha_n)}}\right)^{-1} \leq \epsilon$. For $n \geq n_0$, equations (\ref{eq3}) and (\ref{eq4}) yield

  $$\limsup_{\delta \to 1}  v_{\alpha,\beta, \delta} \geq 1-\epsilon,\; {\rm and}\; \liminf_{\delta \to 1} v_{\alpha,\beta, \delta} \leq \epsilon, $$
    \noindent and the proof of proposition \ref{pro5} is complete.

     \newpage
     
 \subsection{A zero-sum hidden stochastic game} \label{sub4}
 
The MDP and games   considered so far have perfect information and infinite state space. We now mimic the previous construction with a hidden stochastic game with 6 states and 6 public signals.

$\alpha$ and $\beta$ being parameters in $(0,1)$, the HSG  $\Gamma^*(\alpha, \beta)$ is defined as follows. 
The set of states is $K=\{(1,1),(1,0),(2,1),(2,0),1^*,0^*\}$, and the set of public signals is $S=\{s_1, s'_1, s_1^*, s_2, s'_2, s_0^*\}.$ The players perfectly observe past actions and public signals, but not current states.  As in the previous stochastic game, the sum of the payoffs of the players is constantly 1, and the states $0^*$ and $1^*$ are absorbing. The payoffs only depend on the states, player 1 has payoff  0 in states $0^*$, $(1,0)$ and $(1,1)$,  and payoff 1 in  states $1^*$, $(2,0)$ and $(2,1)$. Each player has  2 actions corresponding to  Wait   and  Jump, action sets are $I=\{W_1,J_1\}$ and $J=\{W_2, J_2\}$.  The  initial probability $\pi$ selects with probability $1$ the state $(2,1)$ and the signal $s_2$, so the players know that at period 1 the game is in state $(2,1)$. Once in the absorbing state $0^*$, resp. $1^*$, the play stays there forever and the public signal is $s_0^*$, resp. $s_1^*$. Transitions from states  $(1,0)$ and $(1,1)$ only depend on the action of player 1,  whereas transitions from $(2,0)$ and $(2,1)$ only depend on the action of player 2, and when we write transitions we will omit the action of the player without   influence. More precisely:

\begin{quote}

\item If player 1 Jumps in state $(1,1)$, the play goes to the absorbing state $0^*$ and the public signal is $s_0^*$, i.e. $q((1,1), J_1)$ selects $(0^*, s_0^*)$ a.s. 

\item If player 1 Jumps in state $(1,0)$, the play goes to   state $(2,1)$ and the public signal is $s_2$, i.e. $q((1,0), J_1)$ selects $((2,1), s_2)$ a.s. 

\item Transition when player 1 Waits in state $(1,1)$: $q((1,1), W_1)$ selects $((1,1),s_1)$ with probability $1-\alpha$, $((1,1),s'_1)$ with probability $\alpha ^2$ and $((1,0),s'_1)$ with probability $\alpha (1-\alpha)$. 

\item Transition when player 1 Waits in state $(1,0)$:  $q((1,0), W_1)$ selects $((1,1),s_1)$ with probability $1-\alpha$, and   $((1,0),s'_1)$ with probability $\alpha$. 
\end{quote}

\noindent Transitions from the states controlled by player 2 are  defined symmetrically: $q((2,1), J_2)$ selects $(1^*, s_1^*)$ a.s.,  $q((2,0), J_2)$ selects $((1,1), s_1)$ a.s.,    $q((2,1), W_2)$ selects $((2,1),s_2)$ with probability $1-\beta$, $((2,1),s'_2)$ with probability $\beta ^2$ and $((2,0),s'_2)$ with probability $\beta (1-\beta)$, and finally $q((2,0), W_2)$ selects $((2,1),s_2)$ with probability $1-\beta$  and   $((2,0),s'_2)$ with probability $\beta$. 

Payoffs are discounted with discount factor $\delta\in [0,1)$.

\begin{center}
\begin{figure}[!h]

\begin{center}
\begin{tikzpicture}[>=stealth',shorten >=1pt,auto,node distance=4cm,thick,main
 node/.style={circle,draw,font=\Large\bfseries}]
 \draw (0,8) node[below]{\text{Player 1}};
 \draw (9,8) node[below]{\text{Player 2}};

\node [draw,text width=0.8cm,text centered,circle] (O) at (0,0) {$(1,1)$};
\node [draw,text width=0.8cm,text centered,circle] (A) at (0,5) {$(1,0)$};
\node [draw,text width=0.8cm,text centered,circle] (B) at (0,-3) {$0^*$};

\node [draw,text width=0.8cm,text centered,circle] (O') at (9,0) {$(2,1)$};
\node [draw,text width=0.8cm,text centered,circle] (A') at (9,5) {$(2,0)$};
\node [draw,text width=0.8cm,text centered,circle] (B') at (9,-3) {$1^*$};

\draw[->,>=latex] (O) to[loop right] node[midway, right] {$W_1,1-\alpha,s_1$}(O);
\draw[->,>=latex] (O) to[bend right] node[midway, right] {$W_1,\alpha(1-\alpha),s_1'$}(A);
\draw[->,>=latex] (O) edge[loop left] node[midway, left] {$W_1,\alpha^2,s_1'$}(O);
\draw[->,>=latex] (A) edge[loop above] node[midway, above] {$W_1,\alpha,s_1'$}(A);
\draw[->,>=latex] (A) to[bend right] node[midway, left] {$W_1,1-\alpha,s_1$} (O);

\draw[->,>=latex] (A) to[] node[very near start, above right] {$J_1,s_2$}(O');
\draw[->,>=latex] (O) to[] node[midway, right] {$J_1,s_0^*$}(B);

\draw[->,>=latex] (O') to[loop right] node[midway, right] {$W_2,1-\beta,s_2$}(O');
\draw[->,>=latex] (O') to[bend right] node[midway, right] {$W_2,\beta(1-\beta),s_2'$}(A');
\draw[->,>=latex] (O') edge[loop left] node[midway, left] {$W_2,\beta^2,s_2'$}(O');
\draw[->,>=latex] (A') edge[loop above] node[midway, above] {$W_2,\beta,s_2'$}(A');
\draw[->,>=latex] (A') to[bend right] node[midway, left] {$W_2,1-\beta,s_2$} (O');

\draw[->,>=latex] (A') to[] node[very near start, above left] {$J_2,s_1$}(O);
\draw[->,>=latex] (O') to[] node[midway, right] {$J_2,s_1^*$}(B');

\end{tikzpicture}

\caption{Transitions in $\Gamma^*_{\alpha,\beta}$}
\end{center}
\end{figure}
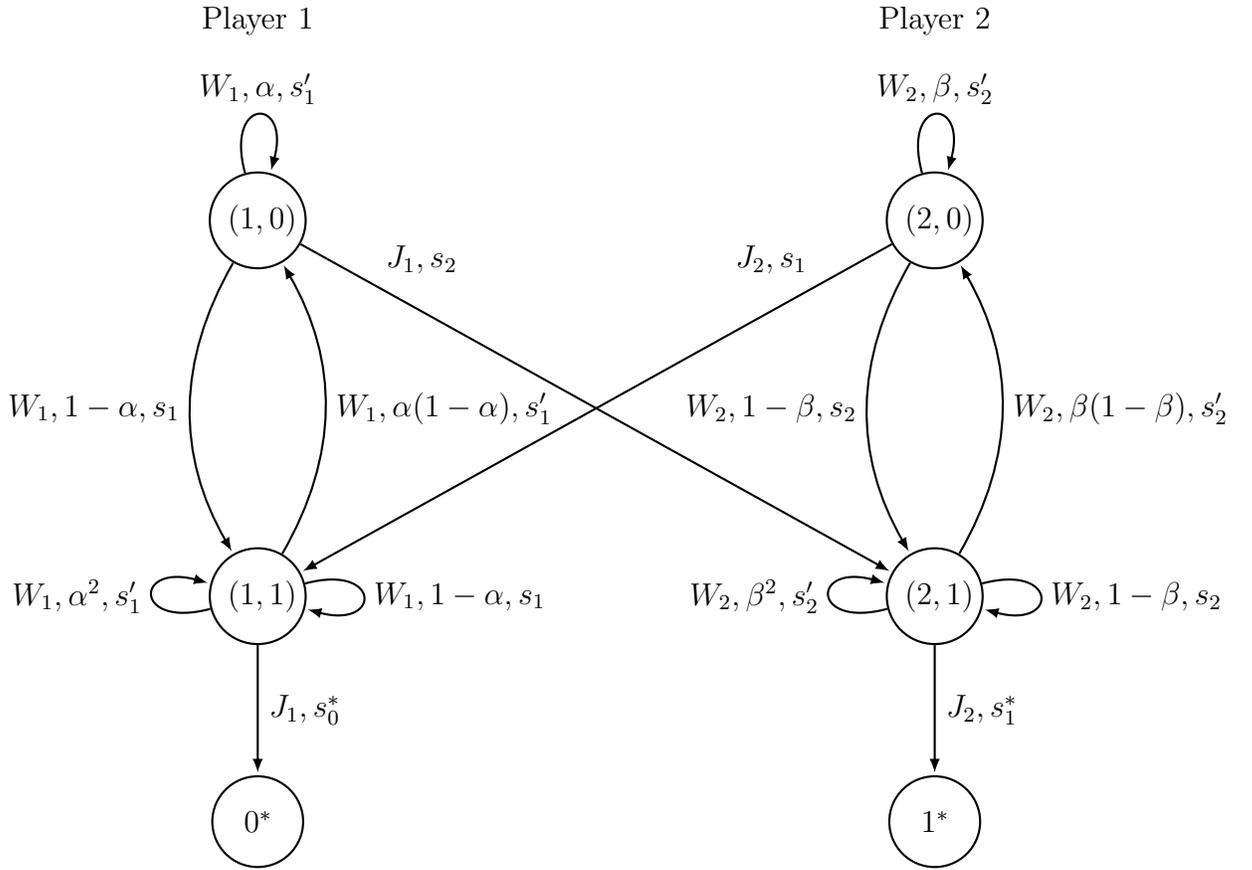
\end{center}

\newpage

Signals in states $(1,0)$ and $(1,1)$ are either  $s_1$ or $s'_1$, and signals  in states $(2,0)$ and $(2,1)$ are either  $s_2$ or $s'_2$. So the public signal always informs the players the element of the partition $\{ \{(1,0), (1,1)\}, \{2,0),(2,1)\}, \{0^*\}, \{1^*\}\}$ that contains the current state, and the game has known payoffs. 

In $\Gamma^*(\alpha, \beta)$, player 1 would like to Jump in state $(1,0)$, and to Wait in state $(1,1)$ but the current state is not fully known to  the players.  Because of the previous partition, the   belief of the players over the current state has at most 2 points  in its support. Suppose this belief corresponds to the state being $(1,1)$ with probability $q$ and   $(1,0)$ with probability $1-q$. The current payoff will be 0, and the transition only depends on player 1's action:
\begin{quote}
\item If player 1 Jumps, the new state is $0^*$ with probability $q$ and $(2,1)$ with probability $1-q$. 

\item If player 1 Waits:   with probability $1- \alpha$ the public signal will be  $s_1$  and by Bayes' rule  the players can deduce that the new state is almost surely $(1,1)$.  With     probability $\alpha$ the public signal is $s'_1$, the probability that the transition selects $((1,1), s'_1)$ is $q \alpha ^2$ so  by Bayes' rule the belief of the players over the new state is :  $(1,1)$ with probability $q \alpha$ and $(1,0)$ with probability $1-q\alpha$. 
\end{quote}

\noindent Consequently the transitions and the payoffs here perfectly mimic those of the stochastic game of subsection \ref{sub3}. The  equivalent stochastic game associated to the HSG $\Gamma^*(\alpha, \beta)$ (see the beginning of section \ref{secHSG}) corresponds to the game $\Gamma(\alpha,\beta)$, up to the addition of the observation of the public signal at the beginning of each period.  This addition plays no role on the payoffs and could only be used as a correlation device for the players, but in a zero-sum context this  has no  influence on the value. We obtain:

\begin{pro} \label{pro6} The value of the $\delta$-discounted  hidden stochastic game $\Gamma^*(\alpha, \beta)$ is  the value $v_{\alpha,\beta, \delta}$ of the $\delta$-discounted stochastic game $\Gamma(\alpha, \beta)$.
\end{pro}

\subsection{A final example} \label{sub5}

Fix $\varepsilon\in (0,\frac{5}{12}]$ and $r$ in $(0, \varepsilon/5)$, we finally construct a non zero-sum HSG $\Gamma$ satisfying the conditions of theorem \ref{thm1}. 
By proposition \ref{pro5}, it is possible to fix  $\alpha$ and $\beta$ such that:  $$\liminf_{\delta\to 1}v_{\alpha,\beta, \delta}< \varepsilon-5r \; {\rm and}\; \limsup_{\delta\to 1}v_{\alpha,\beta, \delta}> {\varepsilon +5r}$$
\noindent And we define:  $\Delta_1=\{\delta \in [1-2r,1), \; v_{\alpha,\beta, \delta}< \varepsilon-5r \}$ and $\Delta_2=\{\delta \in [\frac{1}{1+2r},1), \;  v_{\alpha,\beta, \delta}> {\varepsilon +5r}\}.$ 

%Both $\Delta_1$ anf $\Delta_2$ contain discount factors arbitrarily close to 1. 

Because we want all payoffs of $\Gamma$ to be in $[0,1]$, we   first modify the zero-sum HSG $\Gamma^*(\alpha, \beta)$ of subsection \ref{sub4} by transforming all  payoffs $(1,0)$ into $(1-r,r)$ and all payoffs $(0,1)$ into $(r,1-r)$. That is, we apply the affine increasing transformation $(x\mapsto r+ (1-2r)x)$ to the payoffs, and the game remains constant-sum. We obtain a new HSG $\Gamma_1$ with each payoff in $[r,1-r]$, and the $\delta$-discounted value of this new game is simply $v_\delta=r+ (1-2r)v_{\alpha,\beta, \delta}$. We also define the HSG $\Gamma_2$ as the game $\Gamma_1$ where the identity of the players are exchanged: player 1 in $\Gamma_2$ plays the role of player 2 in $\Gamma_1$, and vice-versa. Plainly, the value of $\Gamma_2$ is $1-v_\delta$.\\

We now define our final HSG  $\Gamma$. 
The states are the 6 states  $(1,1),(1,0),(2,1),$ $(2,0),1^*,0^*$ of $\Gamma_1$, 4 more  states\footnote{There is no need to duplicate states $0^*$ and $1^*$.}  corresponding to the states $(1,1),(1,0),(2,1),$ $(2,0)$ of   $\Gamma_2$, plus 3 extra states $k_1$,  $(\varepsilon, \varepsilon)^*$ and $(1-\varepsilon, 1-\varepsilon)^*$:  $k_1$ is the initial state and is known to the players,  and 
 $(\varepsilon, \varepsilon)^*$ and $(1-\varepsilon, 1-\varepsilon)^*$ are absorbing states where the payoffs will partly depend on the actions played.  
 Actions sets are $I=\{W_1,J_1\}\times \{T,B\}$ and $J=\{ W_2, J_2\}\times \{L,R\}$.  $\Gamma$ is defined as the ``independent sum" of two different games played in parallel, the first game evolving according to the first coordinate of the actions, and the second game evolving according to the second coordinate of the actions. \\
 
1)  At the first period,  the actions of the players determine, through their first coordinate\footnote{At   period 1,  $W_1,   W_2,  J_1, J_2$ should not be interpreted as Wait or Jump.} a continuation game to be played:

$$
\begin {tabular}{cccc}
& \multicolumn{1}{c}{$$} &\multicolumn{1}{c}{$W_2$} &\multicolumn{1}{c}{$J_2$} \\
\cline{3-4}
& \multicolumn{1}{c}{$W_1$} &\multicolumn{1}{|c}{ $(\varepsilon, \varepsilon)^*$ }&
\multicolumn{1}{|c|}{$ \Gamma_2$} \\
\cline{3-4}
& \multicolumn{1}{c}{$J_1$} &\multicolumn{1}{|c}{ $\Gamma_1$} &
\multicolumn{1}{|c|}{$(1-\varepsilon, 1-\varepsilon)^*$}\\
\cline{3-4}
\end{tabular}
$$

\noindent     If $(W_1, W_2)$, resp. $(J_1, J_2)$  is played in period 1, the game reaches the absorbing state $ (\varepsilon, \varepsilon)^*$, resp. $(1-\varepsilon, 1-\varepsilon)^*$. If $(W_1,J_2)$, resp. $(J_1,W_2)$,  is played in period 1, then from period 2 on the hidden stochastic game $\Gamma_2$, resp. $\Gamma_1$,  is played. The payoffs of the first game in period 1 are respectively  defined as $(\varepsilon, \varepsilon)$, $(1-\varepsilon, 1-\varepsilon)$, $(0,0)$ and $(0,0)$  if $(W_1,W_2)$, $(J_1,J_2)$, $(W_1,J_2)$ and $(J_1,W_2)$  is played.  \\

2) In addition, at every period   of $\Gamma$ the players play, through the second coordinate of their actions,  the following bimatrix game $G$, independently of everything else.

$$
\begin {tabular}{cccc}
& \multicolumn{1}{c}{$$} &\multicolumn{1}{c}{$L$} &\multicolumn{1}{c}{$R$} \\
\cline{3-4}
& \multicolumn{1}{c}{$T$} &\multicolumn{1}{|c}{$r,r$} &
\multicolumn{1}{|c|}{$-r,r$} \\
\cline{3-4}
& \multicolumn{1}{c}{$B$} &\multicolumn{1}{|c}{$r,-r$} &
\multicolumn{1}{|c|}{$-r,-r$}\\
\cline{3-4}
\end{tabular}
$$

\vspace{0,5cm}

 In each period,  the payoffs in $\Gamma$ are the sum of the payoffs of the two games. For instance if the state is  $ (\varepsilon, \varepsilon)^*$ and the second components of the actions are $(B,L)$, then the stage payoffs are $\varepsilon +r$ for player 1 and $\varepsilon -r$ for player 2. If at the first period $(J_1,W_2)$ is played then at any subsequent stage the payoffs of the players are the payoffs in $\Gamma_1$ plus the payoffs in $G$. One can easily check that all payoffs lie in $[0,1]$. 
 
 Past actions are perfectly observed. The public signals are those of $\Gamma_1$ or $\Gamma_2$ when these games are played, and we add one  specific public signal for the initial state and each absorbing state $ (\varepsilon, \varepsilon)^* $ and $(1-\varepsilon, 1-\varepsilon)^* $, so that $\Gamma$ has 13 public signals and is a hidden stochastic game with known payoffs. Moreover the game is symmetric between the players. \\

First notice that   in $G$, each player chooses the payoff of the other player, hence any profile is a Nash equilibrium, and the   equilibrium payoff set  of $G$ is the square of feasible payoffs  $[-r,r]^2$.   For each    initial probability and discount factor, the modification of $\Gamma$ where the game  $G$ is removed has a sequential equilibrium yielding  some payoff $(x,y)$. Combining independently such equilibrium with any sequential equilibrium of  the repetition of $G$ gives a sequential equilibrium of $\Gamma$. Then the square centered in $(x,y)$  with side $2r$ is included in the set of sequential equilibrium of $\Gamma$ for this initial probability and discount factor. This proves the third  item of theorem \ref{thm1}. \\

From now on, we consider the game $\Gamma$ with initial state $k_1$.  The idea is quite simple:  for $\delta$ in $\Delta_1$, $v_\delta$ will be significantly  smaller than $\varepsilon$  and all equilibria of $\Gamma$ will play $(W_1,W_2)$ in the first period; whereas for $\delta$ in $\Delta_2$, $v_\delta$ will be much greater   than $\varepsilon$  and all equilibria of $\Gamma$ will play $(J_1,J_2)$ in period 1.   

\begin{pro} \label{lem5} $\;$

1) For $\delta$ in $\Delta_1$, $E_\delta=E'_\delta$ is the square $ [\varepsilon-r, \varepsilon +r]^2$, and this is also the set of communication equilibria of the $\delta$-discounted game, as well as the set of  stationary equilibrium payoffs of the associated stochastic game.

2) For $\delta$ in $\Delta_2$, $E_\delta=E'_\delta$ is the square $ [1-\varepsilon-r, 1-\varepsilon +r]^2$, and this is also the set of communication equilibria of the $\delta$-discounted game, as well as the set of  stationary equilibrium payoffs of the associated stochastic game.\end{pro}

\noindent{\bf Proof:}  First consider, for any discount $\delta$,     the subgame induced by $\Gamma$  after $(J_1,W_2)$ has been played in period  1, discounted from period 2 on. By playing optimally in the $\Gamma_1$ component, player 1 can secure a payoff of $v_\delta-r$, whereas player 2 can secure a payoff of $1-v_\delta-r$. Since the sum of the payoffs is not greater than $1+2r$, all equilibrium payoffs of this subgame lie in the set $[v_\delta-r,v_\delta+3 r]\times [1-v_\delta-r,1-v_\delta+3r]$.  Symmetrically,   equilibrium payoffs of  the subgame induced by $\Gamma$  after $(J_1,W_2)$ has been played in period  1,  belong to  the square  $[1-v_\delta-r,1-v_\delta+3r]\times [v_\delta-r,v_\delta+3r]$. \\

1) Fix    a discount factor $\delta$ in $\Delta_1$. We have $v_\delta=r +(1-2r) v_{\alpha, \beta, \delta}$, so $\delta v_\delta<\varepsilon-4r.$ Consider   a Nash equilibrium $(\sigma, \tau)$ of the $\delta$-discounted game $\Gamma$, and denote by $x$, resp. $y$, the probability that $\sigma$ plays $W_1$, resp. $\tau$ plays $W_2$ at stage 1. We will show that $x=y=1$, and first assume for the sake of contradiction that $x<1$. By playing $W_1$ at period 1 and optimally in $\Gamma_2$ afterwards, player 1 can get a payoff not lower than:
$$A:=y (\varepsilon-r)+(1-y) (\delta (1-v_\delta)-r).$$ 
\noindent This should not exceed the payoff obtained against $\tau$ by playing $J_1$ at period 1 and following $\sigma$ afterwards, and this payoff is not greater than $$B:=y(\delta(v_\delta+3r)+(1-\delta)r) + (1-y) (1-\varepsilon +r),$$ \noindent  because if $y>0$ the continuation strategies after $(J_1,W_2)$ should form a Nash equilibrium of the corresponding subgame. 
Because $\delta v_\delta< \varepsilon -2r (1+\delta)$, we obtain that $\varepsilon-r>\delta(v_\delta+3r)+(1-\delta)r$. Because $\delta v_\delta< \varepsilon -4r$ and $\delta\geq 1-2r$, we have $\delta v_\delta<\varepsilon -2r+\delta -1$, and this implies $\delta (1-v_\delta)-r>1-\varepsilon +r$. Consequently, for all values of $y$ in $[0,1]$ we have $A>B$, which is a contradiction. 
 Hence we obtain $x=1$, and by symmetry $y=1$.  All Nash equilibrium of $\Gamma$ play $W_1$ and $W_2$ in period 1, and the set of Nash equilibrium payoffs $E_\delta$ is included in the square $[\varepsilon-r, \varepsilon +r]^2$. The players can combine $(W_1,W_2)$ in period 1 with the repetition of any given mixed Nash equilibrium of $G$, so any point in the square can be achieved at equilibrium, and $E_\delta=[\varepsilon-r, \varepsilon +r]^2$. Considering sequential equilibria, or introducing a correlation device, even with  communication,  would not modify the above proof. And this is the same with stationary equilibria of the associated stochastic game with state variable the belief   on $K$. This proves 1) of the proposition. \\

2) We proceed similarly for  $\delta$ in $\Delta_2$. We have $v_\delta > r +(1-2r)(\varepsilon +5r)>\varepsilon + 4r$, which implies both: 
$\delta v_\delta> \frac{\varepsilon + 4r}{1 +2r}>\varepsilon +2r$, and $\delta v_\delta> \delta (\varepsilon + 4r) \geq \varepsilon + 2r + \delta(1+2r)-1$.  Let $(\sigma, \tau)$ be  a Nash equilibrium  of the $\delta$-discounted game $\Gamma$, and with  $x$, resp. $y$, being the probability that $\sigma$ plays $W_1$, resp. $\tau$ plays $W_2$, at period  1. Assume for the sake of contradiction that $x>0$. By playing $W_1$ at period 1 and following $\sigma$ afterwards, the payoff of player 1 against $\tau$ is at most:
$$A':= y (\varepsilon+r)+(1-y) ((1-\delta)r+ \delta(1-v_\delta+3r)).$$
\noindent This should not be lower than the payoff obtained by   playing $J_1$ at period 1 and optimally in $\Gamma_1$ afterwards, so not lower than: 
$$B':=  y((1- \delta)(-r) + \delta(v_\delta-r)) + (1-y) (1-\varepsilon -r).$$ 
\noindent Since $\delta v_\delta>  \varepsilon +2r$ and $\delta v_\delta>    \varepsilon + 2r + \delta(1+2r)-1$, we get $B'>A'$, hence a contradiction. We deduce $x=0$, and by symmetry $y=0$. And point  2) of the proposition follows.\\

 Since $\varepsilon+r < 1-\varepsilon -r$,  proposition \ref{lem5} clearly  implies  that no  converging selection of $(E_\delta)_\delta$ exists. \\

\vspace{0,5cm}

We now consider perturbations of the payoffs. Let, for  $\eta\in [0, \frac{r(\varepsilon-5r)}{4})$,   $\Gamma(\eta)$ be a HSG obtained from $\Gamma$ by perturbing each  payoff  by at most $\eta$, and denote by $E_\delta(\eta)$, resp. $E'_\delta(\eta)$,   the corresponding set  of $\delta$-discounted Nash, resp. sequential   equilibrium payoffs with  initial state $k_1$.

 \begin{pro} \label{lem7} $\;$
 
\noindent  1) For all $\delta$ in $\Delta_1$, $ E_\delta(\eta)\subset [\varepsilon-r-2 \eta, \varepsilon +r+ 2 \eta]^2$. $${\it Moreover,}\;\;\lim_{\eta \to 0} \; \; \lim_{\delta \to 1, \delta \in \Delta_1}E'_\delta(\eta) =E_1, \; {\rm and}\; \lim_{\delta \to 1, \delta \in \Delta_1}  \limsup_{\eta \to 0} \; d(E_\delta(\eta),E_1)=0.$$

\noindent 2) For all  $\delta$ in $\Delta_2$, 
$   E_\delta(\eta)\subset [1-\varepsilon-r-2 \eta, 1-\varepsilon +r+ 2 \eta]^2$.  $${\it Moreover,}\;\;\lim_{\eta \to 0} \; \; \lim_{\delta \to 2, \delta \in \Delta_1}E'_\delta(\eta) =E_2, \; {\rm and}\; \lim_{\delta \to 1, \delta \in \Delta_2}  \limsup_{\eta \to 0} \; d(E_\delta(\eta),E_2)=0.$$

\noindent 3) There is no converging selection   $(x_\delta)_\delta$ of $(E_\delta(\eta))_\delta$.  \\

\noindent 4) The game  $\Gamma(\eta)$ has no uniform equilibrium payoff.
\end{pro}

The proof is in the Appendix, and concludes the proof of Theorem \ref{thm1}.

  \section{Appendix}\label{secAppendix}
  
\begin{defi}
Let $A$ be a  subset of the Euclidean space ${\R}^N$. $A$ is \textit{semi-algebraic} if it is defined by a finite number of polynomial inequalities, i.e.  if $A$ is a finite union of sets, each of these sets  being defined as the conjunction of finitely many weak or strict polynomial inequalities. 
\end{defi}

We believe that the following proposition, which is a direct  consequence of the Main Theorem  in \cite{KPV14}, can be useful in several  contexts.
 
\begin{pro}
Let   $(W_{\delta})_{\delta \in [0,1)}$ be a family of non-empty compact subsets of the Euclidean space $\R^N$. Assume that $\left\{(\delta,x)|, \delta\in [0,1), x \in W_{\delta} \right\}$ is a semi-algebraic subset of $\R^{N+1}$. Then when $\delta$ goes to one, $W_{\delta}$ converges for the Hausdorff distance to   a non-empty compact subset $W$ of $\R^N$. \end{pro}
%\begin{proof}
%Applying the Main Theorem in \cite{KPV14} to $T=[0,1)$ and $A=W$, we deduce that $(W_{\delta})$ converges for the Hausdorff metric when $\delta$ goes to $1$.
%%In particular, $(E''_{\delta})$ converges for the Hausdorff metric.
%\end{proof}
  \noindent   \textbf{Proof of Proposition \ref{pro1}}
 Let $W$ be the set of $(\delta,x,y,r) \in [0,1) \times \left(\R^I \times \R^J \times \R^2 \right)^K$ such that $(x,y)$ is a stationary equilibrium in $\Gamma_{\delta}$, and $r$ is the associated payoff equilibrium. Then $(\delta,x,y,r) \in W$ if and only if for all $(k,i,j) \in K \times I \times J$, it satisfies the following inequalities and equalities :
\begin{eqnarray*}
 &\displaystyle \sum_{i' \in I}& x^{i'}(k)=1, \quad x^{i}(k) \geq 0, \quad \sum_{j' \in J} x^{j'}(k)=1, \quad y^{j}(k) \geq 0,
\\
 &\displaystyle \sum_{i' \in I}& x^{i'}(k)\left((1-\delta) u_2(k,i',j)+\delta \sum_{k' \in K}
q^{k'}(k,i',j)r_2(k')\right) \leq r_2(k),
\\
 & \displaystyle \sum_{j' \in J}& y^{j'}(k)\left((1-\delta) u_1(k,i,j')+\delta \sum_{k' \in K}
q^{k'}(k,i,j')r_1(k')\right) \leq r_1(k),
\\
&\displaystyle \sum_{i' \in I}& \sum_{j' \in J} x^{i'}(k) y^{j'}(k) \left((1-\delta) u_2(k,i',j')+\delta \sum_{k' \in K}
q^{k'}(k,i',j')r_2(k')\right) = r_2(k), 
\\
&\displaystyle \sum_{i' \in I}& \sum_{j' \in J} x^{i'}(k) y^{j'}(k) \left((1-\delta) u_1(k,i',j')+\delta \sum_{k' \in K}
q^{k'}(k,i',j')r_1(k')\right) = r_1(k).
\end{eqnarray*}  
Thus $W$ is a semi-algebraic set.

For $\delta \in [0,1)$, let $W_{\delta}:=\left\{(x,y,r) \in \left(\R^I \times \R^J \times \R^2 \right)^K \ | \ (\delta,x,y,r) \in W \right\}$. Then $W_{\delta}$ is non-empty and compact. Applying the preceding proposition, we deduce that $(W_{\delta})$ converges for the Hausdorff metric when $\delta$ goes to $1$.
In particular, $(E''_{\delta})$ converges for the Hausdorff metric.\\

\noindent   \textbf{Proof of Proposition \ref{pro3}} 

1)  Define $\hat{a}=\hat{a}(\alpha, \delta)$ as the integer part of $a^*=a^*(\alpha, \delta)$, we have $v_{\alpha, \delta}\geq s_{\alpha, \delta}(\hat{a}(\alpha, \delta))$.
Since $\hat{a}> a^*-1$, we have $\alpha^{\hat{a}}\leq \sqrt{\frac{1-\delta}{1-\alpha}} \frac{1}{\alpha}\xrightarrow[\delta\to 1]{}0$. Since $\hat{a}\leq a^*$, we have $(1-\delta) {(\alpha \delta)}^{-\hat{a}}\leq  \sqrt{\frac{1-\delta}{1-\alpha}}\,  \delta^{-a^*}\xrightarrow[\delta\to 1]{}0.$ Consequently, $\lim_{\delta \to 1} s_{\alpha, \delta}(\hat{a}(\alpha, \delta))=1$, which implies that $\lim_{\delta \to 1} v_{\alpha, \delta}=1$.\\

We now turn to the proof of conditions 2) and 3) of proposition \ref{pro3}, and start with a lemma. 
  \begin{lem}\label{lem3} For all $\alpha$ and $\delta$ in $(0,1)$,

%1)  For all $\eta>0$, there exists $\delta_0=\delta_0(\alpha)<1$ such that for all $\delta\geq\delta_0$:
%$$s_{\alpha, \delta}(a^*(\alpha,\delta))\geq 1 - (2-\eta) \sqrt{\frac{1-\delta}{1-\alpha}}.$$

\begin{equation}\label{eq1} 1 - 2 \, \delta^{-a^*-1}\,  \sqrt{\frac{1-\delta}{1-\alpha}}\leq s_{\alpha, \delta}(a^*) \leq  1 - 2    \sqrt{\frac{1-\delta}{1-\alpha}} + 3 \frac{1-\delta}{1-\alpha},\end{equation}

%\noindent and for $a\geq 0$ such that $|a-a^*|\geq 1/2$,
\begin{equation}\label{eq2}s_{\alpha,\delta}(a)\leq 1-  \frac{1}{\sqrt{\alpha^{1/2}}}\sqrt{\frac{1-\delta}{1-\alpha}}+ \frac{1-\delta }{1-\alpha} (\alpha + \frac{1}{ \alpha^{1/2}}).\end{equation}\end{lem}

\indent{\bf Proof of lemma \ref{lem3}:} We use  $1-\alpha \delta\geq 1-\alpha$ in the fist line below and $ \delta^{-a^*-1}\geq 1$ in the third line below  to obtain the LHS of (\ref{eq1}):
\begin{eqnarray*}
s_{\alpha, \delta} (a^*)& \geq  &\frac{1-\alpha^{a^*}}{1+\frac{1-\delta}{1-\alpha} \alpha^{-a^*} \delta^{-a^*-1}},\\
& =& \frac{1-\sqrt{\frac{1-\delta}{1-\alpha}}}{1+\sqrt{\frac{1-\delta}{1-\alpha}} \, \delta^{-a^*-1}},\\
 & \geq & \frac{1-\sqrt{\frac{1-\delta}{1-\alpha}}\, \delta^{-a^*-1}}{1+\sqrt{\frac{1-\delta}{1-\alpha}} \, \delta^{-a^*-1}},\\
 & \geq & 1- 2 \sqrt{\frac{1-\delta}{1-\alpha}} \delta^{-a^*-1}.
\end{eqnarray*}

For   inequality (\ref{eq2}), we introduce  $l_{\alpha,\delta}(a)= \frac{1-\alpha^a}{1+\frac{1-\delta}{1-\alpha} \alpha^{-a}\delta^{-a-1}}$.
 If $a\leq a^*-1/4$, we have $\alpha^a\geq \alpha^{a^*-1/4}= \sqrt{\frac{1-\delta}{\alpha^{1/2}(1-\alpha)}}$, and $l_{\alpha,\delta}(a)\leq 1-\alpha^a\leq 1- \sqrt{\frac{1-\delta}{\alpha^{1/2}(1-\alpha)}}$. If $a\geq a^*+1/4$, we have $\alpha^{-a}\geq \alpha^{-a^*-1/4}$ and we write:
\begin{eqnarray*}
l_{\alpha,\delta}(a) & \leq & \frac{1}{1+\frac{1-\delta}{1-\alpha} \alpha^{-a}}\\
 & \leq &  \frac{1}{1+ \sqrt{\frac{1-\delta}{\alpha^{1/2}(1-\alpha)}  }}\\
& \leq & 1- \sqrt{\frac{1-\delta}{\alpha^{1/2}(1-\alpha)}  } + \frac{1-\delta}{\alpha^{1/2}(1-\alpha)}.
\end{eqnarray*}
\noindent And  inequality  (\ref{eq2}) is obtained after noticing that:
\begin{eqnarray*}
s_{\alpha,\delta}(a) & =  & l_{\alpha,\delta}(a) +(1-\delta) \frac{\alpha}{1-\alpha} l_{\alpha,\delta}(a)\\
& \leq &  l_{\alpha,\delta}(a) +(1-\delta) \frac{\alpha}{1-\alpha}.
\end{eqnarray*}
We conclude with the RHS of  (\ref{eq1}). We use  $ \delta^{-a^*-1}\geq 1$ in the first inequality below, and $\frac{1-x}{1+x}\leq 1- 2x + 2x^2$ for all $x\geq 0$ in the third  inequality below:
\begin{eqnarray*}
l_{\alpha,\delta}(a^*)& \leq & \frac{1-\alpha^{a^*}}{1+\frac{1-\delta}{1-\alpha} \alpha^{-a^*}},\\
 &= &  \frac{1-\sqrt{\frac{1-\delta}{1-\alpha}}}{1+\sqrt{\frac{1-\delta}{1-\alpha}}},\\
 & \leq & 1 - 2\sqrt{\frac{1-\delta}{1-\alpha}} + 2 \frac{1-\delta}{1-\alpha}.
 \end{eqnarray*}
 \noindent $s_{\alpha,\delta}(a^*)\leq  l_{\alpha,\delta}(a^*) + \frac{1-\delta}{1-\alpha}$   finally gives the RHS of (\ref{eq1}). This concludes the proof of lemma \ref{lem3}.\\

We now prove point (2) of proposition \ref{pro3}. Fix $\alpha<1/4$, we have $\frac{1}{\sqrt{\alpha}}>2$ so 
for $\delta$ close enough to 1,
$$2 \delta^{-a^*-1}+ \sqrt{\frac{1-\delta}{1-\alpha}}(\alpha +1/\alpha) <  \frac{1}{\sqrt{\alpha}},$$
which implies that: $$1 - 2 \, \delta^{-a^*-1}\,  \sqrt{\frac{1-\delta}{1-\alpha}}>1-  \frac{1}{\sqrt{\alpha}}\sqrt{\frac{1-\delta}{1-\alpha}}+ \frac{1-\delta }{1-\alpha} (\alpha + 1/ \alpha).$$
For   $\delta\in \Delta_1(\alpha)$,  the $a^*$-strategy is available in the MDP, and the previous inequality shows that it is an optimal strategy. $v_{\alpha,\delta}=s_{\alpha, \delta}(a^*)$, and (\ref{eq1}) of lemma \ref{lem3} implies $\lim_{\delta \to 1, \delta \in \Delta_1(\alpha)}\; \frac{1-v_{\alpha,\delta}}{2 \sqrt{\frac{1-\delta}{1-\alpha}}}=1.$

We finally  prove point (3) of proposition \ref{pro3}, and consider $\delta\in \Delta_2(\alpha)$. The pure stationary strategies available in the MDP are $a$-strategies, with $|a-a^*|\geq 1/4$. Point (\ref{eq2}) of lemma \ref{lem3} then implies that:
$v_{\alpha, \delta}\leq 1-  \frac{1}{\sqrt{\alpha}}\sqrt{\frac{1-\delta}{1-\alpha}}+ \frac{1-\delta }{1-\alpha} (\alpha + 1/ \alpha^{1/2})$, hence the result.\\

\noindent   \textbf{Proof of Proposition  \ref{lem7}} 

For any discount factor, the perturbed game  issued from $\Gamma_1$  may no longer be zero-sum, but the quantity  that   player 1 can guarantee (whatever the strategy of the other player) in this game  is close to    $v_\delta$.   More precisely,  in the subgame induced by $\Gamma(\eta)$  after $(J_1,W_2)$ has been played in period  1,    player 1 can secure a payoff of $v_\delta-r-\eta$, whereas player 2 can secure a payoff of $1-v_\delta-r-\eta$. Since the sum of the payoffs is now not greater than $1+2r+2\eta$, all equilibrium payoffs of this subgame lie in the set $[v_\delta-r-  \eta,v_\delta+3 r + 3 \eta]\times [1-v_\delta-r-  \eta,1-v_\delta+3r+3 \eta]$. Symmetrically,  all equilibrium payoffs of the  subgame induced by $\Gamma(\eta)$  after $(W_1,J_2)$ has been played in period  1,   are in the set $[1-v_\delta-r-  \eta,1-v_\delta+3r+3 \eta]\times [v_\delta-r-  \eta,v_\delta+3 r +3 \eta]$.\\

  \noindent $1)$  Fix $\delta$ in $\Delta_1$, we have $v_\delta<r+(1-2r)(\varepsilon-5r)$ and $\delta\geq 1-2r$. This implies:
\begin{equation} \label{eq7} v_\delta\leq \min\{\varepsilon -4 (r+  \eta), \varepsilon-2(r+ \eta) + \delta- 1\}.\end{equation}
\noindent Mimicking the proof of 1) of proposition \ref{lem5}, we obtain $A(\eta)=y (\varepsilon-r-  \eta)+(1-y) (\delta (1-v_\delta)-r-  \eta)$ and $B(\eta)=y(\delta(v_\delta+3r+3\eta)+(1-\delta)(r+ \eta)) + (1-y) (1-\varepsilon +r+  \eta)$, so that $A(\eta)$ and $B(\eta)$ are obtained from the quantities $A$ and $B$ of that lemma by replacing the payoff $r$ by the payoff $r+ \eta$. By inequality (\ref{eq7}), we have $A(\eta)>B(\eta)$. This implies that any $\delta$-discounted Nash equilibrium of $\Gamma(\eta)$ plays $W_1$ and $W_2$ at the first period, and $ E_\delta(\eta)\subset [\varepsilon-r-  \eta, \varepsilon +r+   \eta]^2 $.\\

   Fix now $\eta$ in $(0, \frac{r(\varepsilon-5r)}{2})$. Define $\Gamma(\eta)(W_1,W_2)$ as the subgame obtained  from $\Gamma(\eta)$ after $(W_1,W_2)$ has been played in period 1. $\Gamma(\eta)(W_1,W_2)$ is a repeated game, with stage payoffs  $\eta$-close to the
   bimatrix:
   
 $$
\begin {tabular}{cccccc}
& \multicolumn{1}{c}{$$} &\multicolumn{1}{c}{$(W_2,L)$} &\multicolumn{1}{c}{$(J_2,L)$} &\multicolumn{1}{c}{$(W_2,R)$} &\multicolumn{1}{c}{$(J_2,R)$} \\
\cline{3-6}
& \multicolumn{1}{c}{$(W_1,T)$} &\multicolumn{1}{|c}{$r+ \varepsilon, r+ \varepsilon$} &
\multicolumn{1}{|c|}{$r+ \varepsilon, r+ \varepsilon$}&\multicolumn{1}{|c|}{$-r+ \varepsilon, r+ \varepsilon$} &\multicolumn{1}{|c|}{$-r+ \varepsilon, r+ \varepsilon$}  \\
\cline{3-6}
& \multicolumn{1}{c}{$(J_1,T)$} &\multicolumn{1}{|c}{$r+ \varepsilon, r+ \varepsilon$} &
\multicolumn{1}{|c|}{$r+ \varepsilon, r+ \varepsilon$}&\multicolumn{1}{|c|}{$-r+ \varepsilon, r+ \varepsilon$} &\multicolumn{1}{|c|}{$-r+ \varepsilon, r+ \varepsilon$} \\
\cline{3-6}
& \multicolumn{1}{c}{$(W_1,B)$} &\multicolumn{1}{|c}{$r+ \varepsilon, -r+ \varepsilon$} &
\multicolumn{1}{|c|}{$r+ \varepsilon, -r+ \varepsilon$}&\multicolumn{1}{|c|}{$-r+ \varepsilon, -r+ \varepsilon$} &\multicolumn{1}{|c|}{$-r+ \varepsilon, -r+ \varepsilon$} \\
\cline{3-6}
& \multicolumn{1}{c}{$(J_1,B)$} &\multicolumn{1}{|c}{$r+ \varepsilon,- r+ \varepsilon$} &
\multicolumn{1}{|c|}{$r+ \varepsilon, -r+ \varepsilon$}&\multicolumn{1}{|c|}{$-r+ \varepsilon, -r+ \varepsilon$} &\multicolumn{1}{|c|}{$-r+ \varepsilon, -r+ \varepsilon$} \\
\cline{3-6}
\end{tabular}
$$\\   
  By the   Folk Theorem of Fudenberg and Maskin  (1986), the set $E'_\delta(\eta)(W_1,W_2)$ of sequential equilibrium payoffs of $\Gamma(\eta)(W_1,W_2)$ converges, when $\delta$ goes to 1, to the set of feasible and individually rational payoffs of this game. And this set now converges, when $\eta$ goes to 0, to the square $E_1=[-r+ \varepsilon,r+ \varepsilon]^2$.  Since all sequential equilibria of $\Gamma(\eta)$ play $(W_1,W_2)$ in period 1,   we obtain $\lim_{\eta \to 0} \; \; \lim_{\delta \to 1, \delta \in \Delta_1}E'_\delta(\eta) =E_1.$\\
   
   Consider now the repetition of the bimatrix game $G$. Fix $\varepsilon'>0$, there exists $\delta'$ such that for all $\delta \geq \delta'$ and  any payoff $u$ in $[-r,r]^2$, there exists a  periodic sequence $(i_t,j_t)_t$ of pure action profiles in $\{T,B\}\times \{L,R\}$ such that for all $t_0$, playing the sequence $(i_t,j_t)_{t\geq t_0}$  yields a  $\delta$-discounted payoff $\varepsilon'$-close to $u$. Assume $u=(u_1,u_2)\in [-r+2 \varepsilon',r]^2$ and $\eta<\min \{\varepsilon', \frac{r(\varepsilon-5r)}{2}\}$, we have $u_l-\varepsilon'>-r + \eta$ for each player $l=1,2$. For $\delta\in \Delta_1$, $\delta\geq \delta'$, the   strategy profile where:   $(W_1,W_2)$ is played at stage 1, and for the second component of the actions, the above sequence of pure actions is played, with deviations   punished  by   repeating forever $(J_1,J_2)$, is a Nash equilibrium of the $\delta$-discounted game $\Gamma(\eta)$. Hence   $E_\delta(\eta)$ contains a point $\varepsilon'$-close to $u$, and $d(E_\delta(\eta),E_1)\leq 2 \varepsilon'$. So $\limsup_{\eta \to 0} \; d(E_\delta(\eta),E_1)\leq 2 \varepsilon'$, and $\lim_{\delta \to 1, \delta \in \Delta_1}  \limsup_{\eta \to 0} \; d(E_\delta(\eta),E_1)=0.$\\

\noindent   $2)$  For $\delta$ in $\Delta_2$, we have $\delta v_\delta>\delta(r+(1-2r)(\varepsilon+5r))$. Since $\eta<\frac{1}{2}(1-\varepsilon-5r)$, we have $r+(1-2r)(\varepsilon+5r)> \varepsilon +4( r+ \eta),$ and since $\varepsilon-1 + 2(r +  \eta)<0$, it implies $r+(1-2r)(\varepsilon+5r)> \frac{1}{\delta}(\varepsilon-1 + 2(r +  \eta))+1 +2(r +  \eta).$ So:
 \begin{equation}\label{eq81} \delta v_\delta> \varepsilon-1 + 2 (r+  \eta) + \delta (1+2(r+  \eta)).
 \end{equation}
      
      \noindent Since $\delta\geq\frac{1}{1+2r}$, the above also implies: \begin{equation}\label{eq82}  \delta v_\delta > \varepsilon + 2 (r +  \eta). \end{equation}
      
      We mimick the proof of 2) of proposition  \ref{lem5} and obtain quantities $A'(\eta)= y (\varepsilon+r+  \eta)+(1-y) ((1-\delta)(r+ \eta)+ \delta(1-v_\delta+3r+3 \eta))$, and $B'(\eta)=  y((1- \delta)(-r- \eta) + \delta(v_\delta-r- \eta)) + (1-y) (1-\varepsilon -r-  \eta).$ And the inequalities (\ref{eq81}) and (\ref{eq82}) imply that $B'(\eta)>A'(\eta)$, hence  any $\delta$-discounted Nash equilibrium of $\Gamma(\eta)$ plays $J_1$ and $J_2$ at the first period. 
  The rest of the proof of 2)  is similar to the proof of 1).    \\
  
  \noindent 3) We have  $\varepsilon +r+   \eta< \varepsilon + r (1+  \frac{1}{2}{ \varepsilon-5r})< 1/2$ since  $r<\varepsilon/5$ and $\varepsilon<5/12$.  Hence  there is no converging selection   $(x_\delta)_\delta$ of $(E_\delta(\eta))_\delta$. \\
  
  \noindent 4) It remains to prove that $\Gamma(\eta)$ has no equilibrium payoff, i.e.   that for $\varepsilon'$ small enough, there is no strategy profile which is an $\varepsilon'$-equilibrium of all discounted games $\Gamma(\eta)$ with high enough discount factors.
  
 We proceed by contradiction, and assume that for each   $\varepsilon'>0$, on can find  a discount $\delta_{\varepsilon'}$ in (0,1), and a strategy profile $(\sigma, \tau)=(\sigma_{\varepsilon'}, \tau_{\varepsilon'})$ which is an $\varepsilon'$-equilibrium of each game  $\Gamma(\eta)$ with discount $\delta>\delta_{\varepsilon'}$.  Denote by $x=x_{\varepsilon'}$, resp. $y=y_{\varepsilon'}$, the probability that $\sigma$ plays $W_1$, resp. $\tau$ plays $W_2$ at stage 1. The $\delta$-discounted payoff of player 1 induced by $(\sigma_\varepsilon, \tau_\varepsilon)$ is by definition:$$g_1^\delta(\sigma, \tau)=\E_{\sigma, \tau}\left((1-\delta) \sum_{t=1}^\infty \delta^{t-1} u_1(k_t,i_t,j_t)\right).$$
  \noindent We denote  by $g_1^\delta(\sigma, \tau|W_1,W_2)$ the conditional payoff of player 1 given that $(W_1,W_2)$ is played at period 1, that is:   
      $$E_{\sigma, \tau} \left. \left((1-\delta) \sum_{t=1}^\infty \delta^{t-1} u_1(k_t,i_t,j_t) \right| (i_1=(W_1,T)\; {\rm or} \; (W_1,B))\; {\rm and}\; (j_1=(W_2,L)\; {\rm or}\; (W_2,R)) \right).$$
\noindent And we similarly define $g_1^\delta(\sigma, \tau|W_1,J_2)$, $g_1^\delta(\sigma, \tau|J_1,W_2)$, $g_1^\delta(\sigma, \tau|J_1,J_2)$ and  similar quantities for player 2's payoff. We have:
\begin{eqnarray*}
g_1^\delta(\sigma, \tau) = &xy g_1^\delta(\sigma, \tau|W_1,W_2)+ x(1-y) g_1^\delta(\sigma, \tau|W_1,J_2)\\
& + (1-x)y g_1^\delta(\sigma, \tau|J_1,W_2)+(1-x)(1-y) g_1^\delta(\sigma, \tau|J_1,J_2). 
\end{eqnarray*}
     Because player 1 can secure the  payoff $v_\delta$ in the game $\Gamma_1$, the fact that  $(\sigma, \tau)$ is an $\varepsilon'$-equilibrium implies that:  $$g_1^\delta(\sigma, \tau|J_1,W_2)\geq \delta v_\delta -(r+  \eta) -\frac{\varepsilon'}{(1-x)y}.$$ Similarly,
     $g_1^\delta(\sigma, \tau|W_1,J_2)\geq \delta (1-v_\delta) -(r+  \eta) -\frac{\varepsilon'}{x(1-y)}$,  
     $g_2^\delta(\sigma, \tau|W_1,J_2)\geq \delta  v_\delta -(r+  \eta) -\frac{\varepsilon'}{x(1-y)}$, and $g_2^\delta(\sigma, \tau|J_1,W_2)\geq \delta  (1- v_\delta) -(r+  \eta) -\frac{\varepsilon'}{(1-x)y}$. Since $g_1^\delta(\sigma, \tau|W_1,J_2)+g_2^\delta(\sigma, \tau|W_1,J_2)\leq 1+2r+2 \eta$, we obtain: 
    \begin{eqnarray}
    g_1^\delta(\sigma, \tau|W_1,J_2)& \leq & 1+3 (r +   \eta)  - \delta  v_\delta   +\frac{\varepsilon'}{x(1-y)}\label{eq8}\\
     g_1^\delta(\sigma, \tau|J_1,W_2)& \leq & 1+3 (r +  \eta) - \delta  (1-v_\delta) +  \frac{\varepsilon'}{y(1-x)}\label{eq9}\end{eqnarray}
     
a)      By definition $(\sigma, \tau)$ is an $\varepsilon'$-equilibrium, so playing $J_1$ at period 1 then optimally afterwards  against $\tau$ should not increase player1's payoff by more than  $\varepsilon'$, i.e;
     $$y \sup_{\sigma'}    g_1^\delta(\sigma', \tau|J_1,W_2)+ (1-y)\sup_{\sigma'}    g_1^\delta(\sigma', \tau|J_1,J_2) \leq \varepsilon'+ g_1^\delta(\sigma, \tau).$$
      This implies:
      $$xy  \sup_{\sigma'}    g_1^\delta(\sigma', \tau|J_1,W_2)+x(1-y) \sup_{\sigma'}    g_1^\delta(\sigma', \tau|J_1,J_2)\leq \varepsilon' + xy g_1^\delta(\sigma, \tau|W_1,W_2)+ x(1-y)g_1^\delta(\sigma, \tau|W_1,J_2).$$
     We have  $g_1^\delta(\sigma, \tau|W_1,W_2)\leq \varepsilon +r+  \eta$, $ \sup_{\sigma'}    g_1^\delta(\sigma', \tau|J_1,W_2)\geq \delta v_\delta -r- \eta$ and  $\sup_{\sigma'}    g_1^\delta(\sigma', \tau|J_1,J_2)\geq 1-\varepsilon-r- \eta$. Together with  inequality (\ref{eq8}), it implies: 
 \[    xy (\delta v_\delta -r-2\eta) + x(1-y) (1-\varepsilon -r-  \eta) \leq 2 \varepsilon'+ xy (\varepsilon +r+  \eta) + x(1-y) ( 1 + 3 (r+  \eta) -\delta v_\delta).\]
      \noindent Rearranging terms, the above equation is  equivalent to:
      $$2\varepsilon' + 2x (r+   \eta)(2-y) \geq x (\delta v_\delta-\varepsilon).$$
      
 $x=x_{\varepsilon'}$ and $y=y_{\varepsilon'}$ depend on $\varepsilon'$.     Consider $\delta$ in $\Delta_2$, we have $v_\delta > \varepsilon +4( r+  \eta)$. So there exists $\varepsilon''>0$, independent from $\varepsilon'$, such that for all $\delta$ high enough in $\Delta_2$:
 $$2\varepsilon' + 2x_{\varepsilon'} (r+   \eta)(2-y_{\varepsilon'}) \geq 4x_{\varepsilon'} (r+   \eta) + x_{\varepsilon'}\varepsilon''.$$
      Passing to the limit  gives: $$x_{\varepsilon'} \xrightarrow[\varepsilon' \to 0]{}0.$$
      And by symmetry between the players, we also have $\lim_{\varepsilon'\to 0} y_{\varepsilon'}=0$.\\
      
      b) We finally write that playing $W_1$ at period 1 then optimally afterwards  against $\tau$ should not increase player 1's payoff by more than  $\varepsilon'$, i.e;
     $$y \sup_{\sigma'}    g_1^\delta(\sigma', \tau|W_1,W_2)+ (1-y)\sup_{\sigma'}    g_1^\delta(\sigma', \tau|W_1,J_2) \leq \varepsilon'+ g_1^\delta(\sigma, \tau).$$
 This implies:
   $y(1-x)  \sup_{\sigma'}    g_1^\delta(\sigma', \tau|W_1,W_2)+ (1-y)(1-x) \sup_{\sigma'}    g_1^\delta(\sigma', \tau|W_1,J_2) $
      
       $\leq \varepsilon' + (1-x)y g_1^\delta(\sigma, \tau|J_1,W_2)+ (1-x)(1-y)g_1^\delta(\sigma, \tau|J_1,J_2).$\\
   
   We have  $g_1^\delta(\sigma, \tau|J_1,J_2)\leq 1-\varepsilon +r+  \eta$, $ \sup_{\sigma'}    g_1^\delta(\sigma', \tau|W_1,W_2)\geq \varepsilon-r- \eta$ and  $\sup_{\sigma'}    g_1^\delta(\sigma', \tau|W_1,J_2)\geq \delta(1-v_{\delta})-r- \eta$. Together with  inequality (\ref{eq9}), the above  implies : 
   $$2 \varepsilon'+ 2 (r+   \eta) (1-x)(1+y) \geq (1-x) (\varepsilon-1 + \delta(1-v_\delta)).$$
      
For $\delta \in \Delta_1$, we have $v_\delta< \varepsilon-4 (r+  \eta)$ so  for all $\delta$ high enough in $\Delta_1$:  $\varepsilon-1 + \delta(1-v_\delta)\geq 4 (r+  \eta)$    and we obtain: $$\frac{\varepsilon'}{r+  \eta}+(1-x)(1+y) \geq 2(1-x).$$
\noindent We finally get a contradiction since $\lim_{\varepsilon'\to 0} x=\lim_{\varepsilon'\to 0} y=0$.


\begin{thebibliography}{99}


  
       \bibitem{APS90}
Abreu D., D.~Pearce et  E.~Stacchetti.
\newblock Toward a theory of discounted repeated games with imperfect
monitoring.
  {\em Econometrica}, 58, 1041--1063, 1990.
  
\bibitem{AumSha94}
Aumann R.J.  and L.~S. Shapley.
\newblock Long-term competition---{A} game theoretic analysis.
\newblock In N.~Megiddo, editor, {\em Essays on game theory}, pages 1--15.
  Springer-Verlag, New-York, 1994.


  
%  \bibitem{BK85}
%Benoit J-P.~ and V.~Krishna  
%\newblock Finitely repeated games,
%\newblock {\em Econometrica}, {\bf 53}, 905--922, 1985.

%\bibitem{BK87}
%Benoit J.-P.~ and V.~Krishna.
%\newblock Nash equilibria of finitely repeated games,
%\newblock {\em International Journal of Game Theory}, {\bf 16}, 197--204, 1987.

\bibitem{BK76}
  Bewley T.  and E. Kohlberg. \newblock The asymptotic theory of stochastic games. \newblock{\em Mathematics of Operations Research}, {\bf 1}:197208, 1976.
 
 
 \bibitem{BGV13}
 Bolte J., S. Gaubert and G. Vigeral. \newblock Definable zero-sum stochastic games. Arxiv:1301.1967, 2013. 
 
 \bibitem{Dutta}
  Dutta P.K. \newblock A Folk Theorem for stochastic games. \newblock{\em Journal of Economlc Theory}, {\bf 66}:1--32, 1995.



\bibitem{Fo85}
Forges F.
\newblock An Approach to Communication Equilibria.    {\em Econometrica}, 54,
1375--1385, 1986.

%\bibitem{ForgesMertensNeymanCounterExampleFolk1986}
%F.~Forges, J.-F. Mertens, and A.~Neyman.
%\newblock A counterexample to the folk theorem with discounting.
%\newblock {\em Economics Letters}, 20:7--7, 1986.

  
%  \bibitem{FMN86} Forges F., J.-F. Mertens and  A.~Neyman (1986)
%\newblock A counterexample to the Folk theorem with discounting,
%\newblock {\em Economic Letters}, {\bf  20}, 7.

\bibitem{FM86} Fudenberg D. and E.~Maskin.
 The Folk Theorem in repeated games with discounting or with
incomplete information.
\newblock {\em Econometrica}, {\bf 54}, 533--554, 1986.

\bibitem{FL91} Fudenberg D. and D.~Levine.
 An approximate Folk Theorem with imperfect private information.
\newblock {\em Journal of Economic Theory}, {\bf 54}, 26--47, 1991.

\bibitem{FY11} Fudenberg D. and Y.~Yamamoto.
The folk theorem for irreducible stochastic games with imperfect public monitoring.
\newblock {\em Journal of Economic Theory}, {\bf 146}, 1664--1683, 2011.

%\bibitem{FudLev94}
%Fudenberg D. and D. Levine.
%\newblock Efficiency and observability with long-run and short-run players.
%\newblock {\em Journal of Economic Theory}, 62:103--135, 1994.


\bibitem{FudLevMas94}
Fudenberg D., D. Levine, and E.~Maskin.
\newblock The folk theorem with imperfect public information.
\newblock {\em Econometrica}, 62:997--1039, 1994.

\bibitem{FLT04}
Fudenberg D., D. Levine, and S.~Takahashi.
\newblock Perfect public equilibrium when players are patient.
\newblock {\em Games and Economic Behavior}, 61:27--49, 2007.

\bibitem{PG83}
 Georghiou C., Philippou  A.N and Philippou G.N.
 \newblock  A Generalized Geometric Distribution
and some of its Properties. {\it Statistics \& Probability Letters}, 1, 171-175, 1983. 

\bibitem{G14}
 Gimbert, H., Renault, J., Sorin, S., Venel, X. and Zielonka, W.
\newblock
  On the values of repeated games with signals 
arXiv preprint arXiv:1406.4248, 2014.

%\bibitem{Gos95}
%Gossner O. 
%\newblock The Folk Theorem for finitely repeated games with mixed strategies.
%\newblock {\em International Journal of Game Theory}, {\bf 24}, 95--107, 1995.

%\bibitem{GT07}
% Gossner O. and  T.~Tomala.
%\newblock Secret correlation in repeated games with imperfect monitoring.
%  {\em Mathematics of Operations Research}, 32, 413--424, 2007.
\bibitem{HSTV} H\"{o}rner J.,  Sugaya T.,  Takahashi S.  and Nicolas Vieille. \newblock { Recursive Methods in Discounted Stochastic Games: An Algorithm for $\delta \to1$ and a Folk Theorem}. \newblock {\em  Econometrica},  {\bf 79}, 1277-1318, 2011,

%\bibitem{HM96}
%  Hillas J.  and M. Liu.
%\newblock Repeated Games with Partial Monitoring: the Stochastic Signaling Case.
%\newblock Game Theory and Information 9605001, EconWPA, 1996.

%\bibitem{EHO}
% Ely J., J.~H\"{o}rner, and W.~Olszewski.
%\newblock Belief-free equilibria in repeated games.
%\newblock {\em Econometrica}, 73:377--415, 2005.

%
%\bibitem{KM98}
% Kandori M. and H.~Matsushima.
%\newblock Private observation, communication and collusion.
%\newblock {\em Econometrica}, 66:627--652, 1998.

% 
\bibitem{KPV14}
Kocel-Cynk B., Pawlucki W. and Valette A.
\newblock A short geometric proof that Hausdorff limits are definable in any o-minimal structure. \newblock{\em Advances in Geometry},  14-1,  49--58, 2014.

%\bibitem{Leh89}
% Lehrer E.
%\newblock Lower Equilibrium Payoffs in Two-Player Repeated Games with Non-observable Actions.
%  {\em International Journal of Game Theory}, 18, 57--89, 1989.

 


\bibitem{Leh90}
 Lehrer E.
\newblock Nash equilibria of $n$-player repeated games with semi-standard
information.
  {\em International Journal of Game Theory}, 19, 191--217, 1990.


%\bibitem{Leh92a}
% Lehrer E.
%\newblock Correlated Equilibria in two-Player Repeated Games with
%non-Observable Actions.
%  {\em Mathematics of Operations Research}, 17, 175--199, 1992a.

\bibitem{Leh92b}
 Lehrer E.
\newblock On the Equilibrium Payoffs Set of
two-Player Repeated Games with Imperfect
Monitoring.
  {\em International Journal of Game Theory}, 20, 211--226, 1992a.

\bibitem{Leh92c}
 Lehrer E.
\newblock Two-player repeated games with nonobservable actions and observable payoffs.
  {\em Mathematics of Operations Research}, 17, 200--224, 1992b.


 
%\bibitem{MMS}  Mailath G., S. Matthews and  T. Sekiguchi. \newblock Private strategies in Finitely Repeated Games with Imperfect Public Monitoring. \newblock Contributions to Theoretical Economics, Vol. 2, article 2, 2002.

%\bibitem{MS06}
% Mailath G. and L.~Samuelson.
%\newblock {\em Repeated Games and Reputations: Long-Run Relationships}.
%\newblock Oxford University Press, 2006.


\bibitem {MSZa94}
 Mertens J-F., S.~Sorin et  S.~Zamir.
\newblock Repeated games.
\newblock CORE discussion paper  9420, Louvain-la-Neuve, 1994.


 

\bibitem{mye}
 Myerson R.
\newblock Multistage games with communication.
  {\em Econometrica},
54, 323--358, 1986.


\bibitem{neyman}
 Neyman A. \newblock Real algebraic tools in stochastic games. \newblock {\em Stochastic Games and Applications}. Chapter 6,  NATO Science Series,  A. Neyman and S. Sorin eds, 2003.
 

%\bibitem{}  Renault J. and T. Tomala. \newblock  Repeated Proximity Games.       
%{\it  International Journal of Game Theory},
%2,  539--559, 1998. 

 \bibitem{}  Renault  J. and T. Tomala. \newblock Communication equilibria    in repeated games with imperfect monitoring.  \newblock{\em Games and Economic Behavior}. {\bf 49}, 313--344,  2004. 

\bibitem{}  Renault  J. and T. Tomala. \newblock General Properties of Long-Run Supergames. \newblock{\em Dynamic Games and Application}. {\bf 1}, 319--350, 2011. 


%\bibitem{Rub77}
% Rubinstein, A.
%\newblock Equilibrium in supergames.
%\newblock Center for Research in Mathematical Economics and Game Theory,
%  Research Memorandum 25, 1977.
  
\bibitem{Rub1994} Rubinstein A. \newblock {Equilibrium in supergames", N.Meggido (ed.)}, {\it  Essays in Game Theory in Honor of Michael Maschler}, Springer-Verlag, 17-28, 1994.
  
%    \bibitem{Seki2001}  Sekiguchi T.  \newblock A negative result in finitely repeated games with product structure. {\it Economic Letters},   74,   67-70, 2001.
%  
%  \bibitem{Seki2005}  Sekiguchi T.  \newblock Uniqueness of equilibrium payoffs in finitely repeated games with imperfect monitoring. {\it Japanese Economic Review},   56,   317-331, 2005.
  
\bibitem{Shapley}
Shapley L.S.
\newblock Stochastic games.
\newblock {\em Proceedings of the National Academy of Sciences of the U.S.A.}, 39, 1095-1100, 1953.

%\bibitem{So92}
% Sorin, S.
%\newblock Repeated games with complete information.
%\newblock In Handbook of Game Theory with Economic Applications, R.J. Aumann and S. Hart (eds.), Elsevier Science Publishers, 71-107, 1992.

\bibitem{So86}
 Sorin S.
\newblock On repeated games with  complete information.
\newblock {\em Mathematics of Operations Research}, {\bf 11}, 147--160, 1986.

\bibitem {Sor84}
 Sorin S.
\newblock Asymptotic properties of a non-zero sum stochastic game.  
\newblock {\em   International Journal of Game Theory},
98:296--303, 1984.

\bibitem{Venel}
Venel X. \newblock Commutative stochastic games. Preprint, 2012. \newblock 

\bibitem{Wen94} Wen, Q. \newblock The ``Folk Theorem" for repeated games with Complete Information. \newblock{\em Econometrica}, {\bf 62}, 949--954, 1994.

\bibitem{V00a} Vieille N.   Two-player stochastic games I: a reduction. 
\textit{Isra\"{e}l Journal of Mathematics, } {119}, 55-91, 2000a.

\bibitem{V00b} Vieille N.   Two-player stochastic games II: the case of
recursive games, \textit{Isra\"{e}l Journal of Mathematics, } {119},
93-126, 2000b.



\bibitem{BZ2013} Ziliotto B.  \newblock Zero-sum repeated games: counterexamples to the existence of the asymptotic value and the conjecture Maxmin= Lim $v_N$. arXiv:1305.4778, 2013.
\end{thebibliography}
 \end{document}